\documentclass[12pt]{amsart}
\usepackage{latexsym,amssymb,amsrefs,graphicx}
\usepackage[all]{xy}

\newtheorem{lemma}{Lemma}

\newtheorem{proposition}[lemma]{Proposition}
\newtheorem{theorem}[lemma]{Theorem}
\newtheorem{corollary}[lemma]{Corollary}
\newtheorem{definition}[lemma]{Definition}

\theoremstyle{definition}
\newtheorem{example}[lemma]{Example}

\begin{document}

\title[Morse-Bott Homology]{Morse-Bott Homology}

\author{Augustin Banyaga}
\address{Department of Mathematics \\
         Penn State University \\
         University Park, PA 16802}
\email{banyaga@math.psu.edu}

\author{David E. Hurtubise}
\address{Department of Mathematics and Statistics\\
         Penn State Altoona\\
         Altoona, PA 16601-3760}
\email{Hurtubise@psu.edu}


\begin{abstract}
We give a new proof of the Morse Homology Theorem by constructing
a chain complex associated to a Morse-Bott-Smale function that reduces 
to the Morse-Smale-Witten chain complex when the function is Morse-Smale
and to the chain complex of smooth singular $N$-cube chains when the 
function is constant. We show that the homology of the chain complex 
is independent of the Morse-Bott-Smale function by using compactified moduli 
spaces of time dependent gradient flow lines to prove a Floer-type
continuation theorem.
\end{abstract}

\maketitle


\section{Introduction}

\subsection{Overview}
Let $Cr(f) = \{p \in M |\, df_p = 0 \}$ denote the set of critical points of a 
smooth function $f:M \rightarrow \mathbb{R}$ on a  smooth $m$-dimensional manifold
$M$. A critical point $p \in Cr(f)$ is said to be nondegenerate if and only if 
$df$ is transverse to the zero section of $T^\ast M$ at $p$. In local coordinates 
this is equivalent to the condition that the $m \times m$ Hessian matrix
$
\left(\frac{\partial^2 f}{\partial x_i \partial x_j}\right)
$
has rank $m$ at $p$.  If all the critical points of $f$ are non-degenerate, then 
$f$ is called a Morse function.

A Morse function $f:M \rightarrow \mathbb{R}$ on a finite dimensional compact 
smooth Riemannian manifold $(M,g)$ is called Morse-Smale if all its stable and 
unstable manifolds intersect transversally. Such a function gives rise to a 
chain complex $(C_\ast(f),\partial_\ast)$, called the Morse-Smale-Witten chain
complex. The chains of this complex are generated by the critical points of $f$,
and the boundary operator is defined by counting gradient flow lines between 
critical points (with sign). The Morse Homology Theorem says that the homology 
of the Morse-Smale-Witten chain complex is isomorphic to the singular homology of
$M$ with integer coefficients \cite{BanLec}, \cite{SchMor}.

\vskip .2in

If $Cr(f)$ is a disjoint union of finitely many connected submanifolds 
and the Hessian of $f$ is non-degenerate in the direction normal 
to $Cr(f)$, then $f$ is called a Morse-Bott function \cite{BotNon}.
In this paper we construct a new chain complex associated to a Morse-Bott
function meeting certain transversality requirements (see Definition \ref{MBStransversality})
that generalizes the Morse-Smale-Witten chain complex. This singular Morse-Bott-Smale chain 
complex has the following form:
$$
\xymatrix{
\ddots & \vdots\\
\cdots & C_1(B_2) \ar@{}[u]|{\oplus} \ar[r]^{\partial_0} \ar[dr]^{\partial_1} \ar[ddr]|(.33){\partial_2} |!{[d];[dr]}\hole & C_0(B_2) \ar[r]^{\partial_0} \ar[dr]^{\partial_1} \ar[ddr]|(.33){\partial_2} |!{[d];[dr]}\hole & 0 & & \\
\cdots & C_2(B_1)\ar@{}[u]|{\oplus} \ar[r]^(.55){\partial_0} \ar[dr]^{\partial_1} & C_1(B_1)\ar@{}[u]|{\oplus} \ar[r]^(.55){\partial_0} \ar[dr]^{\partial_1} & C_0(B_1) \ar@{}[u]|{\oplus} \ar[r]^{\partial_0} \ar[dr]^{\partial_1} & 0 & \\
\cdots & C_3(B_0) \ar@{}[u]|{\oplus} \ar[r]^{\partial_0} & C_2(B_0) \ar@{}[u]|{\oplus} \ar[r]^{\partial_0} & C_1(B_0) \ar@{}[u]|{\oplus} \ar[r]^{\partial_0} & C_0(B_0) \ar@{}[u]|{\oplus}  
 \ar[r]^-{\partial_0} & 0\\
\cdots &  {C}_3(f) \ar@{}[u]|{\|} \ar[r]^{\partial} & {C}_2(f) \ar@{}[u]|{\|} \ar[r]^{\partial} & {C}_1(f) \ar@{}[u]|{\|} \ar[r]^{\partial} & {C}_0(f) \ar@{}[u]|{\|} \ar[r]^-{\partial} & 0
}
$$
where $C_p(B_i)$ is the group of ``$p$-dimensional chains'' in the critical 
submanifolds of Morse-Bott index $i$, and the boundary operator is defined 
as a sum of homomorphisms $\partial = \partial_0 \oplus \cdots \oplus \partial_m$. 

When the function $f:M \rightarrow \mathbb{R}$ is Morse-Smale the critical
set $B_i$ is a discrete set of points for all $i=0,\ldots,m$, and the
groups $C_p(B_i)$ are trivial for all $p > 0$. When the function is constant
the entire manifold $M$ is a critical submanifold of Morse-Bott index zero.
In this case $B_i = \emptyset$ for all $i > 0$, and the groups $C_p(B_i)$ 
are trivial for all $i > 0$. These two cases appear in the above 
diagram of a general Morse-Bott-Smale chain complex as follows.
$$
\xymatrix{
&  \ddots \ar[dr]^{\partial_1}& & & & \\
&  &  C_0(B_2) \ar[dr]^{\partial_1} &  & & \\
&  & & C_0(B_1) \ar[dr]^{\partial_1} &  & \\
\cdots \ar[r]^{\partial_0} & C_3(B_0) \ar[r]^{\partial_0} & C_2(B_0) 
\ar[r]^{\partial_0} & C_1(B_0)  \ar[r]^{\partial_0} & C_0(B_0) \ar[dr]^{\partial_1}
\ar[r]^-{\partial_0} & 0\\
& & & & & 0
}
$$
In the first case the homomorphism $\partial_1$ is the Morse-Smale-Witten
boundary operator, and in the second case the homomorphism $\partial_0$ 
is comparable to the boundary operator on singular cubical chains found in
\cite{MasABa}. Thus, the Morse-Bott-Smale chain complex provides a means 
of interpolating between the Morse-Smale-Witten chain complex and
the chain complex of singular cubical chains. Moreover, a Floer-type continuation
theorem shows that the homology of the Morse-Bott-Smale chain complex
is independent of the Morse-Bott-Smale function $f:M \rightarrow \mathbb{R}$
and gives a new proof of the Morse Homology Theorem.

\vskip .2in

Heuristically, the homomorphisms $\partial_j:C_p(B_i) \rightarrow C_{p+j-1}(B_{i-j})$
for $j=1,\ldots ,m$ can be described as follows. Imagine that every 
topological space under consideration has a prefered triangulation.
Then a map $\sigma:P \rightarrow B_i$ from some triangulated space
$P$ of dimension $p$ would determine a chain in the singular
chain group $S_p(B_i;\mathbb{Z})$. Using $\sigma:P \rightarrow B_i$
we could form a fibered product with respect to the beginning point map 
$\partial_-:\overline{\mathcal{M}}(B_i,B_{i-j}) \rightarrow B_i$
from the compactified moduli space of gradient flow lines from $B_i$ to $B_{i-j}$,
and then by projecting onto the second component of the fibered product 
and composing with the endpoint map $\partial_+:\overline{\mathcal{M}}(B_i,B_{i-j})
\rightarrow B_{i-j}$ we would have a map 
$$
\partial_j:P \times_{B_i} \overline{\mathcal{M}}(B_i,B_{i-j}) 
\stackrel{\pi_2}{\longrightarrow} \overline{\mathcal{M}}(B_i,B_{i-j})
\stackrel{\partial_+}{\longrightarrow} B_{i-j}.
$$
The dimension of $\overline{\mathcal{M}}(B_i,B_{i-j})$ is
$b_i + j - 1$, where $b_i = \text{dim }B_i$, and (assuming some
transversality conditions) the dimension of the fibered
product $P \times_{B_i} \overline{\mathcal{M}}(B_i,B_{i-j})$
is $p+j-1$.  Thus, if we were to assume that the fibered product
also has a preferred triangulation, then the above map would determine a 
singular chain in the singular chain group $S_{p+j-1}(B_{i-j};\mathbb{Z})$.

Of course, a topological space does not necessarily have
a prefered triangulation. One approach to making the above
heuristic discussion rigorous would be to show that all the spaces
under consideration can be triangulated in some way and then prove
that the homology of the resulting complex does not depend on
the chosen triangulations. However, this approach would involve
picking triangulations on an uncountably infinite number of spaces,
namely all the fibered products $P \times_{B_i} \overline{\mathcal{M}}(B_i,B_{i-j})$.
This approach seems quite daunting in light of the fact that there are simple
diagrams of polyhedra and piecewise linear maps that are not triangulable:
$$
Y \stackrel{\sigma_1}{\longleftarrow} X \stackrel{\sigma_2}{\longrightarrow} Z.
$$
That is, there may not exist triangulations of $X, Y$, and $Z$ with
respect to which both $\sigma_1$ and $\sigma_2$ are simplicial
\cite{BryTri}. (See also Example \ref{badfibered} in Subsection \ref{manicorners}
of this paper.)

\vskip .2in

In light of these observations, our approach to constructing
a singular Morse-Bott-Smale chain complex is not based on picking
triangulations, but on defining singular chain complexes where
the generators are allowed to have a variety of domains. We generalize 
the usual singular chain complexes based on maps from the standard
$p$-simplex $\Delta^p \subset \mathbb{R}^{p+1}$ (or the unit $p$-cube $I^p$)
to allow maps from the faces of an $N$-cube (for some large $N$),
the endpoint maps from compactified moduli spaces of gradient flow lines,
and the endpoint maps from their fibered products. We call these 
generalizations of $\Delta^p$ abstract topological chains, and we refer
to the corresponding singular chains as singular topological chains.
Thus, the maps
$$
\partial_j:P \times_{B_i} \overline{\mathcal{M}}(B_i,B_{i-j}) 
\stackrel{\pi_2}{\longrightarrow} \overline{\mathcal{M}}(B_i,B_{i-j})
\stackrel{\partial_+}{\longrightarrow} B_{i-j}.
$$
for $j=1,\ldots ,m$ described above are singular topological
chains in our new framework (without any need for a triangulation
on the domain). However, the homomorphism $\partial_0$ 
is no longer the ``usual'' boundary operator on singular chains, 
since the generators have a variety of domains.

Of course, a singular chain complex with maps defined on
a variety of domains will contain more maps than a singular
chain complex where all the $p$-chains are defined on the domain
$\Delta^p$, and some redundancy may appear in the larger
chain complex. Hence, we define degeneracy relations 
to identify singular chains that are ``essentially'' the same
but defined on different domains. For example, fix a large integer 
$N$ and let $C_p$ be the $p$-dimensional faces of the unit $N$-cube $I^N$.
There are a large number of elements in $C_p$ (which are all ``essentially''
the same), and different elements may support singular chains that 
are ``essentially'' the same. For instance, any element $P \in C_p$
is homeomorphic to the standard unit $p$-cube $I^p$ via a homeomorphism
$\alpha:P \rightarrow I^p$, and the usual singular boundary operator
on a singular $p$-cube $\sigma:I^p \rightarrow B$ in a topological 
space $B$ induces a singular boundary operator on the singular $C_p$-space 
$\sigma \circ \alpha:P \rightarrow B$ such that the boundary operators
commute with the homeomorphism $\alpha:P \rightarrow I^p$.  So, 
$\sigma:I^p \rightarrow B$ and $\sigma \circ \alpha:P \rightarrow B$ 
are ``essentially'' the same, even though they are defined on 
different domains.

\vskip .1in

More precisely, to construct our singular Morse-Bott-Smale chain complex we 
begin by defining the set $C_p$ of allowed domains of degree $p$ for all $p$, and then 
we define a boundary operator on the free abelian group generated by the elements
of $C_p$, i.e. $S_p = \mathbb{Z}[C_p]$. This boundary operator
$$
\partial:S_p \rightarrow S_{p-1}
$$
induces a boundary operator
$$
\partial:S_p(B) \rightarrow S_{p-1}(B)
$$
on singular $C_p$-chains in $B$ as follows. If $\sigma:P \rightarrow B$ is a 
singular $C_p$-space in $B$, then $\partial(\sigma)$ is given by the formula
$$
\partial(\sigma) = \sum_k n_k \sigma|_{P_k}
$$ 
where 
$$
\partial(P) = \sum_k n_k P_k
$$
and $n_k = \pm 1$ for all $k$.  The degeneracy relations define 
a subgroup $D_p(B)\subseteq S_p(B)$ of degenerate singular topological 
chains in $B$ such that $\partial(D_p(B)) \subseteq D_{p-1}(B)$, and hence 
the boundary operator on singular $C_p$-chains in $B$ descends to a boundary
operator on the quotient groups:
$$
\partial:S_p(B)/D_p(B) \rightarrow S_{p-1}(B)/D_{p-1}(B).
$$

This describes our framework for defining the homomorphism $\partial_0$ 
in the Morse-Bott-Smale chain complex. In fact, in the above diagram for a general 
Morse-Bott-Smale chain complex we have
$$
C_p(B_i) = S_p^\infty(B_i)/D^\infty_p(B_i)
$$
where $S_p^\infty(B_i) \subset S_p(B_i)$, $D_p^\infty(B_i) 
\subset D_p(B_i)$, and $\partial_0:C_p(B_i) \rightarrow C_{p-1}(B_i)$ 
is $(-1)^{p+i}$ times the boundary operator on singular topological
chains developed in Section \ref{topologicalchains} of this paper (see
Definitions \ref{MBdegree}, \ref{boundary}, and \ref{MBScomplex}). The 
homomorphisms $\partial_1,\ldots ,\partial_m$ are then defined as above using 
fibered products of compactified moduli spaces, and it is shown that they preserve
the groups of degenerate singular topological chains (see Lemma \ref{boundarydegenerate}). 
Hence there are induced homomorphisms
$$
\partial_j:S_p^\infty(B_i)/D_p^\infty(B_i) \rightarrow 
S_{p+j-1}^\infty(B_{i-j})/D_{p+j-1}^\infty(B_{i-j})
$$
for all $j=0,\ldots , m$, and these induced homomorphisms fit
together to form the boundary operator $\partial = \partial_0 \oplus \cdots
\oplus \partial_m$ in the Morse-Bott-Smale chain complex.

\vskip .1in


\subsection{Comparison with other approaches}
The Morse-Bott-Smale chain complex constructed in this paper is inspired 
by several other authors (\cite{AusMor}, \cite{FraTheA}, \cite{FukFlo}, \cite{LatGra}, 
\cite{LiuOnt}, and \cite{RuaBot}) who have studied both the finite dimensional
case addressed here and the infinite dimensional version of Floer-Bott homology.

\vskip .1in

However, the approach taken in this paper is fundamentally different and, 
in our opinion, more straightforward and more like singular homology than
previous approaches. For instance, other authors use the machinery of spectral
sequences, cf. \cite{AusMor}, \cite{FukFlo}, and \cite{RuaBot}.  In contrast,
the approach taken here does not require spectral sequences (although several 
of the diagrams found in Section \ref{MBcomplex} are inspired by and resemble 
first quadrant spectral sequences). 

Moreover, as a novel feature, the Morse-Bott-Smale chain complex 
defined in Section \ref{MBcomplex} provides a common framework for
encoding both smooth singular cubical chains and Morse-Smale-Witten chains.
In particular, it reduces to the chain complex of smooth singular $N$-cube
chains when the function is constant (Example \ref{constantfunction}), and it reduces to 
the Morse-Smale-Witten chain complex when the function is Morse-Smale 
(Example \ref{MSfunction}). This is quite different from the approaches that
count flow lines with cascades \cite{FraTheA}, where the chain complex 
reduces to the Morse-Smale-Witten chain complex for both Morse-Smale
functions and constant functions.

\vskip .1in

To construct a Morse-Bott-Smale chain (or cochain) complex one must
first choose a category for the objects in the chain (or cochain) 
complex. For instance, Austin and Braam \cite{AusMor} use differential
forms on the critical submanifolds to define the groups $C_p(B_i)$ in 
their Morse-Bott-Smale cochain complex, and they define the homomorphisms
$\partial_1,\ldots ,\partial_m$ in their cochain complex using 
integration along the fibers. They then use the machinery of spectral 
sequences to prove that the cochain complex they construct computes 
the de Rahm cohomology of the manifold. 

In contrast, Latschev \cite{LatGra} views the stable and unstable manifolds 
(and their fibered products) as currents. Latschev shows (following \cite{HarFin})
that these currents determine an operator (from differential forms to currents)
that is chain homotopic to the inclusion. After showing that this operator 
can be extended to singular chains that are ``sufficiently transverse" to 
the unstable manifolds, Latschev shows that the image of the extended
operator is the set of ``stable bundles of smooth chains" in the critical 
submanifolds. A formula for the differential on the chain complex whose 
groups consist of stable bundles over chains in the critical submanifolds 
is given, and it is proved that the homology of this (Morse-Bott-Smale)
chain complex is isomorphic to the singular homology of the manifold 
with integer coefficients.

Another approach to Morse-Bott homology can be found in a paper by Fukaya \cite{FukFlo},
where the Floer-Bott homology of the Chern-Simons functional is studied. The groups
in Fukaya's version of Morse-Bott homology are generated by piecewise smooth
maps from objects he calls ``abstract geometric chains" and attributes to Gromov \cite{GroFil}.
Similarly, Ruan and Tian's version of  ``Bott-type symplectic Floer cohomology" \cite{RuaBot}
is based on these same ``abstract geometric chains", which they attribute to Fukaya 
\cite{FukFlo}. According to Fukaya, an ``abstract geometric chain" is a finite simplicial 
complex (see Definition 1.14 of \cite{FukFlo}) that satisfies certain conditions, including
a codimension 2 boundary condition, reminiscent of the conditions used to define ``pseudocycles"
(see \cite{KahPse}, \cite{McDJho}, and \cite{SchEqu}). The boundary operator in the chain complex
described by Fukaya is defined on the set of ``transversal geometric chains" by taking fibered
products with moduli spaces of gradient flow lines, and Theorem 1.2 of \cite{FukFlo} says that
there is a spectral sequence associated to a ``generic'' Morse-Bott function which converges 
to the singular homology of the manifold.

\vskip .1in

In this paper, we construct our Morse-Bott-Smale chain complex over the
category of compact oriented smooth manifolds with corners. 
This is a natural choice of category for a Morse-Bott-Smale
chain complex because the compactified moduli spaces of gradient flow 
lines of a Morse-Bott-Smale function are compact oriented 
smooth manifolds with corners (see Theorem \ref{compactification}). 
Moreover, we prove in Subsection \ref{manicorners} that all the relevant
fibered products are compact smooth manifolds with corners.  Hence,
the fibered product constructions used in this paper stay within the chosen
category for the objects in the chain complex.

\vskip .1in


\subsection{Outline of the paper}
In Section \ref{MorseChain} we recall some basic facts about the Morse-Smale-Witten 
chain complex, and in Section \ref{transversalitysection} we recall some basic
facts about Morse-Bott functions and state the Morse-Bott-Smale 
transversality condition.

\vskip .1 in

Section \ref{topologicalchains} establishes the framework we use to
define the homomorphism $\partial_0$ in the Morse-Bott-Smale chain 
complex constructed in Section \ref{MBcomplex}.

\vskip .1 in 

In Subsection \ref{generalabstract} we introduce the notion of abstract 
topological chains and singular topological chains, which generalize the 
usual singular chain complexes based on maps from the standard $p$-simplex 
$\Delta^p\subset \mathbb{R}^{p+1}$ or the unit $p$-cube $I^p$. Abstract 
topological chains are introduced in order to provide the flexibility of 
constructing singular chain complexes with maps from a variety of domains. 
To construct the Morse-Bott-Smale chain complex in Section \ref{MBcomplex}
we need to allow maps from the faces of an $N$-cube (for some large $N$),
the compactified moduli spaces of gradient flow lines, and their 
fibered products.

\vskip .1in

The formalism of abstract topological chains and singular topological chains
is very general and provides a natural framework for some of the formulas
found in \cite{FukFlo}.

\vskip .1in

In Subsection \ref{Ncubes} we show how the faces of an $N$-cube can be 
viewed as abstract topological chains. Using the degeneracy relations in
Definition \ref{cubedegenerate}, this produces the chain complex of
singular $N$-cube chains, which is comparable to the singular cubical 
chain complex of \cite{MasABa} (see Theorem \ref{cubehomology}).
The reader should note that in Subsection \ref{Ncubes} the singular topological 
chain group $S_p(B)$ is generated by maps from the $p$-dimensional faces of the unit 
$N$-cube. So, the maps in $S_p(B)$ are defined on several different domains.
However, these domains are identified via the subgroup $D_p(B)$ of degenerate 
singular topological chains. Hence, every equivalence class in $S_p(B)/D_p(B)$ 
contains a representative defined on the standard $p$-cube
$I^p \times (0,\ldots ,0) \subseteq I^N$.

In Subsection \ref{fiberedsub} we define the fibered product of singular
topological chains and show that the fibered product of singular topological
chains is an abstract topological chain. The formal constructions in
Subsection \ref{fiberedsub} are used in Subsection \ref{modulisingular} where 
the abstract topological chain structure on compactified moduli spaces
of gradient flow lines is defined. This abstract topological chain structure 
is defined in terms of fibered products (see Definition \ref{moduliboundary}),
and the general results proved in Subsection \ref{fiberedsub} are used
to prove that the degree and boundary operator for compactified moduli
spaces of gradient flow lines satisfy the axioms for abstract topological
chains (see Lemma \ref{modulicomplex}). 

\vskip .1in 

In Section \ref{MBcomplex} we construct our Morse-Bott-Smale chain 
complex $(C_\ast(f),\partial)$. The construction begins by defining
the set $C_p$ of allowed domains for the singular topological chains. 
The set of allowed domains $C_p$ includes the $p$-dimensional faces 
of the unit $N$-cube $I^N$ and the $p$-dimensional connected components
of fibered products with compactified moduli spaces of gradient flow
lines. Subsection \ref{manicorners} is devoted to the proof of Lemma
\ref{basiscorners}, which says that the elements of $C_p$ are
compact oriented smooth manifolds with corners.

In Subsection \ref{MBchaindef} we define the chain complex 
$(\tilde{C}_\ast(f),\partial)$. The groups in the Morse-Bott-Smale 
chain complex $(C_\ast(f),\partial)$ are defined in Definition \ref{MBScomplex}
as the quotient of the groups in $(\tilde{C}_\ast(f),\partial)$ by the
degenerate singular topological chains. These degenerate singular topological
chains are defined in Subsection \ref{degMB} using the orientation
conventions given in Subsection \ref{orientationssub}.  The degenerate
singular topological chains defined in Subsection \ref{degMB} include the 
degenerate singular $N$-cube chains from Subsection \ref{Ncubes} together
with additional chains having domains that are components of
fibered products. The degenerate singular topological chains are used to
identify chains that are ``essentially'' the same, but defined on different
domains.  

The chain complex $(\tilde{C}_\ast(f),\partial)$ defined in Subsection 
\ref{MBchaindef} has a natural filtration determined by the Morse-Bott
index, and the boundary operator decomposes into a sum of differentials
$\partial = \oplus_{j=0}^m \partial_j$, where $\partial_j$ is a homomorphism that 
decreases the Morse-Bott index by $j$. In fact, in Subsection \ref{MBchaindef} 
the homomorphisms $\partial_j$ are defined first, and then the full boundary 
operator is defined to be $\partial = \oplus_{j=0}^m \partial_j$. Thus, 
the chain complex $(\tilde{C}_\ast(f),\partial)$ looks like a first quadrant 
spectral sequence where the homomorphisms $\partial_j$ are all defined on the
$E^0$ term (see Definition \ref{MBScomplex}). Note however that, in general, 
$\partial_j^2 \neq 0$ for $j>0$ even though $\partial^2 =0$ by
Corollary \ref{Cfischain}.

The homomorphism $\partial_0$ in Subsection \ref{MBchaindef} is defined in
terms of the singular topological chain structure discussed in Section 
\ref{topologicalchains}, and the homomorphisms $\partial_1,\ldots ,\partial_m$
are defined using fibered products of compactified moduli spaces of gradient
flow lines. It is shown in Lemma \ref{boundarydegenerate} that these
homomorphisms preserve the groups of degenerate singular topological
chains, and hence they induce homomorphisms (still denoted by 
$\partial_0,\ldots ,\partial_m$) on the groups in the Morse-Bott-Smale
chain complex $(C_\ast(f),\partial)$.  Thus, Corollary \ref{Cfischain}
implies that $\partial^2 = 0$ for the Morse-Bott-Smale chain complex
$(C_\ast(f),\partial)$. In Subsection \ref{examples} we give several
examples of computing Morse-Bott homology, including the cases where the
function is constant (Example \ref{constantfunction}) and where the function
is Morse-Smale (Example \ref{MSfunction}).

\vskip .1in

Section \ref{continuationsection} is devoted to proving a Floer-type
continuation theorem (Theorem \ref{homologyindependence}) that 
shows that the homology of the Morse-Bott-Smale chain complex 
defined in Section \ref{MBcomplex} is independent of the Morse-Bott-Smale 
function $f:M \rightarrow \mathbb{R}$ used to define the complex. 
The proof follows standard arguments involving compactified moduli 
spaces of time-dependent gradient flow lines that can be found in 
\cite{AusMor}, \cite{BarLag}, \cite{FloCoh}, \cite{FloAnI}, \cite{SchMor},
and \cite{WebThe}. (An outline of the proof of the continuation theorem
can be found at the beginning of Section \ref{continuationsection}.) 
When the Morse-Bott-Smale function $f:M \rightarrow \mathbb{R}$ is constant
the entire manifold $M$ is a critical submanifold of Morse-Bott index zero,
and the Morse-Bott-Smale chain complex defined in Section \ref{MBcomplex} 
reduces to the chain complex of smooth singular $N$-cube chains (Example 
\ref{constantfunction}). Hence, Theorem \ref{homologyindependence} 
implies that, for any Morse-Bott-Smale function $f:M \rightarrow \mathbb{R}$, 
the homology of the Morse-Bott-Smale chain complex associated to $f$
is isomorphic to the singular homology of the manifold with integer
coefficients.  Moreover, when the function $f:M \rightarrow \mathbb{R}$ 
is Morse-Smale the complex reduces to the Morse-Smale-Witten chain complex
(Example \ref{MSfunction}). Thus, the Morse Homology Theorem follows 
as a corollary to the results in this paper.

\vskip .1in

Other authors have shown that the theory of Morse-Bott homology can be
applied to a variety of problems in algebraic topology, Floer homology,
gauge theory, and quantum cohomology.  For instance, applications
to equivariant homology/cohomology and cup products are discussed in 
\cite{AusMor} and \cite{LatGra}, applications to the Floer homology
of the Chern-Simons functional and Donaldson polynomials are
discussed in \cite{AusEqu} and \cite{FukFlo}, and applications to
symplectic Floer homology and quantum cohomology are discussed in
\cite{LiuOnt} and \cite{RuaBot}. Although we do not discuss applications
of the Morse-Bott-Smale chain complex constructed in Sections
\ref{topologicalchains} and \ref{MBcomplex}, the techniques 
developed in this paper were designed with the above applications in mind.
Thus, we expect that our approach to Morse-Bott homology can be applied 
to problems such as the ones listed above and can be readily 
extended to the infinite dimensional case of Floer-Bott homology.


\section{The Morse-Smale-Witten chain complex}\label{MorseChain}

In this section we briefly recall the construction of the Morse-Smale-Witten
chain complex and the Morse Homology Theorem.  For more details we
refer the reader to \cite{BanLec}.

Let $f:M \rightarrow \mathbb{R}$ be a Morse function on a finite dimensional 
compact smooth Riemannian manifold $(M,g)$. The \textbf{stable manifold} of
a critical point $p \in Cr(f)$ is defined to be
$$
W^s(p) = \{ x\in M | \lim_{t \rightarrow \infty} \varphi_t(x) = p \}
$$
where $\varphi_t$ is the 1-parameter group of diffeomorphisms generated by
minus the gradient vector field, i.e. $-\nabla f$.  Similarly, the 
\textbf{unstable manifold} of $p$ is defined to be
$$
W^u(p) = \{ x\in M | \lim_{t \rightarrow -\infty} \varphi_t(x) = p \}.
$$
The Stable/Unstable Manifold Theorem for a Morse Function says that the 
tangent space at $p$ splits as 
$$
T_pM = T^s_pM \oplus T_p^uM
$$
where the Hessian is positive definite on $T_p^sM$ and negative definite
on $T_p^uM$.  Moreover, the stable and unstable manifolds of $p$ are the surjective 
images of smooth embeddings
\begin{eqnarray*}
E^s: T_p^sM & \rightarrow & W^s(p) \subseteq M\\
E^u: T_p^uM & \rightarrow & W^u(p) \subseteq M. 
\end{eqnarray*}
Hence, $W^s(p)$ is a smoothly embedded open disk of dimension $m - \lambda_p$, 
and $W^u(p)$ is a smoothly embedded open disk of dimension $\lambda_p$
where $m$ is the dimension of $M$ and $\lambda_p$ is the index of the critical
point $p$, i.e. the dimension of $T^u_pM$.

If the stable and unstable manifolds of a Morse function $f:M \rightarrow
\mathbb{R}$ all intersect transversally, then the function is called
\textbf{Morse-Smale}. The next result describes the abundance of Morse-Smale
functions.

\begin{theorem}[Kupka-Smale Theorem]\label{KupkaSmale}
If $(M,g)$ is a finite dimensional compact smooth Riemannian manifold, then the space of 
all $C^r$ Morse-Smale functions on $M$ is a dense subspace of $C^r(M,\mathbb{R})$ for any 
$2 \leq r \leq \infty$.
\end{theorem}

\noindent
If $f:M \rightarrow \mathbb{R}$ is a Morse-Smale function, then using
Palis' $\lambda$-Lemma one can show that there is a partial ordering
on the critical points of $f$ defined as follows.

\begin{definition}\label{partialorder}
Let $p$ and $q$ be critical points of $f:M \rightarrow \mathbb{R}$. We say 
that $q$ is \textbf{succeeded} by $p$, $q \succeq p$, if and only if 
$W(q,p) = W^u(q) \cap W^s(p) \neq \emptyset$, i.e. there exists a gradient
flow line from $q$ to $p$.
\end{definition}

\noindent
Note that $W(q,p)$ is an embedded submanifold of $M$ of dimension
$\lambda_q  - \lambda_p$.  Hence, $q \succeq p$ implies that
$\lambda_q \geq \lambda_p$.  As other consequences of Palis' 
$\lambda$-Lemma we have the following.

\begin{theorem}\label{W(q,p)closure}
If $p$ and $q$ are critical points of $f:M \rightarrow \mathbb{R}$ such
that $q \succeq p$, then
$$
\overline{W(q,p)} = \overline{W^u(q)} \cap \overline{W^s(p)} \ \  = 
\bigcup_{q \succeq \tilde{q} \succeq \tilde{p} \succeq p} W(\tilde{q},\tilde{p})
$$
where the union is over all critical points between $q$ and $p$.
\end{theorem}

\begin{corollary}\label{rel1compact}
If $p$ and $q$ are critical points of relative index one, i.e.
if $\lambda_q - \lambda_p = 1$, and $q \succeq p$, then
$$
\overline{W(q,p)} = W(q,p) \cup \{p,q\}.
$$
Moreover, $W(q,p)$ has finitely many components, i.e. the number of gradient
flow lines from $q$ to $p$ is finite.
\end{corollary}

If we assume that $M$ is oriented and we choose an orientation for each of the 
unstable manifolds of $f$, then there is an induced orientation on the stable 
manifolds.  The preceding corollary shows that we can then define an integer 
$n(q,p)$ associated to any two critical points $p$ and $q$ of relative index one
by counting the number of gradient flow lines from $q$ to $p$ with signs
determined by the orientations.  
The \textbf{Morse-Smale-Witten chain complex} is defined to be the chain complex 
$(C_\ast(f),\partial_\ast)$ where $C_k(f)$ is the free abelian group
generated by the critical points $q$ of index $k$ and the boundary
operator $\partial_k:C_k(f) \rightarrow C_{k-1}(f)$ is given by
$$
\partial_k(q)\ \ = \sum_{p \in Cr_{k-1}(f)} n(q,p)p.
$$

\begin{theorem}[Morse Homology Theorem]\label{Morsehomology}
The pair $(C_\ast(f),\partial_\ast)$ is a chain complex, and its homology
is isomorphic to the singular homology $H_\ast(M;\mathbb{Z})$.
\end{theorem}

\noindent
Note that the Morse Homology Theorem implies that the homology of
$(C_\ast(f),\partial_\ast)$ is independent of the Morse-Smale
function $f:M \rightarrow \mathbb{R}$, the Riemannian metric $g$,
and the orientations.


\section{Morse-Bott functions and transversality}\label{transversalitysection}

Let $f:M \rightarrow \mathbb{R}$ be a smooth function on a finite 
dimensional compact smooth Riemannian manifold $(M,g)$. Assume 
that the critical point set
$$
Cr(f) = \{p \in M |\, df_p = 0 \} = \coprod_{\alpha=1}^n B_\alpha
$$
is a finite disjoint union of connected submanifolds $B_\alpha$ in $M$,
called \textbf{critical submanifolds}. For any $p \in B$, a critical submanifold,
the tangent space of $M$ splits as
$$
T_p M \approx T_p B \oplus \nu_p(B)
$$
where $\nu_p(B)$ denotes the normal bundle of $B$ in $M$. For
$V,W \in T_p M$, the Hessian is a symmetric bilinear form defined
by $Hess_p(f)(V,W) = V_p(\tilde{W}(f))$ where $\tilde{W}$ is any
extension of $W$ to a neighborhood of $p$. If $V\in T_p B$ then
$Hess_p(f)(V,W) = 0$ because $\tilde{W}_x(f) = 0$ for all $x \in B$.
Thus the Hessian determines a symmetric bilinear form on the normal 
space
$$
Hess_p^\nu(f):\nu_p(B) \times \nu_p(B) \rightarrow \mathbb{R}.
$$

\begin{definition}\label{MorseBottDefinition}
A function $f:M \rightarrow \mathbb{R}$ is said to be \textbf{Morse-Bott} 
if and only if the critical point set $Cr(f)$ is a finite disjoint 
union of connected submanifolds of $M$ and for each connected submanifold 
$B \subseteq Cr(f)$ the normal Hessian $Hess_p^\nu(f)$ is non-degenerate 
for all $p\in B$. The \textbf{index} $\lambda_p$ of a critical point 
$p \in B \subseteq Cr(f)$ is defined to be the maximal dimension of a 
subspace of $\nu_p(B)$ on which $Hess_p^\nu(f)$ is negative definite.
\end{definition}

For a proof of the following lemma see Section 3.5 of \cite{BanLec} or \cite{BanAPr}.

\begin{lemma}[Morse-Bott Lemma] \label{MorseBottLemma}
Let $f:M \rightarrow \mathbb{R}$ be a Morse-Bott function, and let $B$ be a
connected component of the critical set $Cr(f)$.  For any $p \in B$ there
is a local chart of $M$ around $p$ and a local splitting of the normal
bundle of $B$
$$
\nu_\ast(B) = \nu^+_\ast(B) \oplus \nu^-_\ast(B)
$$
identifying a point $x \in M$ in its domain to $(u,v,w)$ where $u \in B$,
$v \in \nu_\ast^+(B)$, $w \in \nu_\ast^-(B)$ such that within this chart
$f$ assumes the form
$$
f(x) = f(u,v,w) = f(B) + |v|^2 - |w|^2.
$$
\end{lemma}

\noindent
Note that the Morse-Bott Lemma shows that if $B$ is connected, then $\lambda_p$
is independent of $p \in B$. Hence, we may also refer to $\lambda_p$ as the
index $\lambda_B$ of the critical submanifold $B$.  Moreover, the lemma shows 
that at a critical point $p\in Cr(f)$ the tangent space splits as
$$
T_pM = T_pB \oplus \nu^+_p(B) \oplus \nu^-_p(B)
$$
where $\lambda_p = \mbox{dim } \nu^-_p(B)$.  If we let $b = \mbox{dim } B$ and
$\lambda_p^\ast = \mbox{dim } \nu^+_p(B)$, then we have the fundamental
relation
\begin{eqnarray*}
m & = & b + \lambda_p^\ast + \lambda_p
\end{eqnarray*}
where $m = \mbox{dim }M$.

Let $\varphi_t:M \rightarrow M$ denote the flow associated to $-\nabla f$, 
i.e. $\varphi_t(x) = \gamma (t)$ where $\gamma'(t) = -(\nabla f)(\gamma(t))$ and 
$\gamma(0) = x$. For $p \in Cr(f)$ the \textit{stable manifold} $W^s(p)$ 
and the \textit{unstable manifold} $W^u(p)$ are defined the same as 
they were for a Morse function:
\begin{eqnarray*}
W^s(p) & = & \{ x \in M | \lim_{t \rightarrow \infty} \varphi_t(x) = p \} \\
W^u(p) & = & \{ x \in M | \lim_{t \rightarrow -\infty} \varphi_t(x) = p \}.
\end{eqnarray*}
However, for a Morse-Bott function we can also consider the stable and unstable
manifolds of a critical submanifold $B\subseteq Cr(f)$.  These are defined
to be
\begin{eqnarray*}
W^s(B) & = & \bigcup_{p \in B} W^s(p)\\
W^u(B) & = & \bigcup_{p \in B} W^u(p).
\end{eqnarray*}
For the following see Proposition 3.2 of \cite{AusMor}.
\begin{theorem}[Stable/Unstable Manifold Theorem for a Morse-Bott Function]\label{stablemanifold}
The stable and unstable manifolds $W^s(B)$ and $W^u(B)$ are the surjective
images of smooth injective immersions $E^+:\nu^+_\ast(B) \rightarrow M$ and 
$E^-:\nu^-_\ast(B) \rightarrow M$.  There are smooth endpoint maps 
$\partial_+:W^s(B) \rightarrow B$ and $\partial_-:W^u(B) \rightarrow B$ 
given by $\partial_+(x) = \lim_{t \rightarrow \infty} \varphi_t(x)$ and 
$\partial_-(x) = \lim_{t \rightarrow -\infty} \varphi_t(x)$ which when 
restricted to a neighborhood of $B$ have the structure of locally trivial 
fiber bundles.
\end{theorem}

\begin{definition}[Morse-Bott-Smale Transversality]\label{MBStransversality}
A Morse-Bott function $f:M \rightarrow \mathbb{R}$ is said to satisfy the 
\textbf{Morse-Bott-Smale transversality} condition with respect to a given 
metric on $M$ if and only if for any two connected critical submanifolds $B$ and $B'$,
$W^u(p)$ intersects $W^s(B')$ transversely, i.e. $W^u(p) \pitchfork W^s(B')$, 
for all $p\in B$.
\end{definition}

\smallskip\noindent
Note: By the Kupka-Smale Theorem (Theorem \ref{KupkaSmale}) we can always
find a Morse-Smale function as close as we like to a given Morse-Bott
function. However, there are situations where a certain Morse-Bott function
has some desired properties that will be lost if it is perturbed. In those
cases it might be preferable to look for a Riemannian metric such that the
function is Morse-Bott-Smale with respect to the chosen metric.  Unfortunately,
it is not always possible to perturb the Riemannian metric to make a given
Morse-Bott function satisfy the Morse-Bott-Smale transversality condition.
See Section 2 of \cite{LatGra} for some interesting counterexamples.

\smallskip
If the gradient flow of $f:M \rightarrow \mathbb{R}$ does satisfy the
Morse-Bott-Smale transversality condition, then the space
$$
W(B,B') = W^u(B) \cap W^s(B')
$$
of points $x \in M$ such that $\partial_-(x) \in B$ and $\partial_+(x) \in B'$
is a submanifold of $M$.  The following lemma is an immediate consequence
of transversality.
\begin{lemma}\label{dim}
Suppose that $B$ is of dimension $b$ and index $\lambda_B$ and that $B'$ is of
dimension $b'$ and index $\lambda_{B'}$. Then we have the following where $m =
\mbox{dim }M$:
\begin{eqnarray*}
\mbox{dim }W^u(B) & = & b + \lambda_B \\
\mbox{dim }W^s(B') & = & b' + \lambda_{B'}^\ast = m - \lambda_{B'}\\
\mbox{dim }W(B,B') & = & \lambda_B - \lambda_{B'} + b \quad (\mbox{if } W(B,B') \neq \emptyset).
\end{eqnarray*}
\end{lemma}
\noindent
Note that the dimension of $W(B,B')$ does not depend on the dimension  of the
critical submanifold $B'$.  This fact will be used when we define the
boundary operator in the Morse-Bott-Smale chain complex.

We end this section with the following lemma which says that
a Morse-Bott-Smale function is \textbf{weakly self-indexing} \cite{AusMor},
i.e. the Morse-Bott index is strictly decreasing along gradient flow lines.

\begin{lemma}\label{weakly}
If $f:M \rightarrow \mathbb{R}$ satisfies the Morse-Bott-Smale
transversality condition, then $W(B,B') = \emptyset$ whenever
$\lambda_B \leq \lambda_{B'}$ and $B \neq B'$.  
\end{lemma}

\smallskip\noindent
Proof:  If $B \neq B'$ and there exists an $x \in W(B,B')$, then $x$ lies on a 
one dimensional gradient flow line that begins at $\partial_-(x) \in B$ and ends at 
$\partial_+(x) \in B'$.  Hence, $\mbox{dim }W(\partial_-(x),B') \geq 1$, and
the Morse-Bott-Smale transversality condition implies that 
$\mbox{dim }W(\partial_-(x),B') = \lambda_B - \lambda_{B'} + 0 \geq 1$.
Therefore, $\lambda_B > \lambda_{B'}$.

\begin{flushright}
$\Box$
\end{flushright}


\section{Topological chains and fibered products}\label{topologicalchains}

In this section we introduce the concept of abstract topological chains and
singular topological chains. Abstract topological chains are generalizations
of the standard $p$-simplex $\Delta^p\subset \mathbb{R}^{p+1}$, and singular 
topological chains are generalizations of singular chains where we consider 
continuous maps from abstract topological chains. Applying these definitions
to the faces of an $N$-cube we construct a chain complex that is comparable to 
the singular cubical chain complex of \cite{MasABa}. We also define the fibered
product of singular topological chains and show that the fibered product of singular 
topological chains is an abstract topological chain. The homomorphism $\partial_0$
in the Morse-Bott-Smale chain complex constructed in Section \ref{MBcomplex} 
is defined in terms of abstract topological chains that are constructed from 
fibered products of compactified moduli spaces of gradient flow lines and the 
faces of an $N$-cube. In this section we define the degrees and boundary operators
for these spaces and show that they satisfy the axioms for abstract
topological chains.

\subsection{Some general abstract definitions}\label{generalabstract}

\smallskip
For each integer $p \geq 0$ fix a set $C_p$ of topological spaces, and
let $S_p$ be the free abelian group generated by the elements of $C_p$, 
i.e. $S_p = \mathbb{Z}[C_p]$. Set $S_p = \{0\}$ if $p<0$ or $C_p = \emptyset$.

\begin{definition}\label{topologicalchainsdef}
A \textbf{boundary operator} on the collection $S_\ast$ of groups $\{S_p\}$ 
is a homomorphism $\partial_p:S_p \rightarrow S_{p-1}$ such that 
\begin{enumerate}
\item For $p \geq 1$ and $P \in C_p\subseteq S_p$, $\partial_p(P) = \sum_k n_k P_k$
      where $n_k = \pm 1$ and $P_k\in C_{p-1}$ is a subspace of $P$ for all $k$.
\item $\partial_{p-1}\circ \partial_p:S_p \rightarrow S_{p-2}$ is zero.
\end{enumerate}
We will call $(S_\ast,\partial_\ast)$ a \textbf{chain complex of abstract 
topological chains}. Elements of $S_p$ are called \textbf{abstract topological chains} 
of \textbf{degree} $p$.
\end{definition}

\begin{definition}\label{singulartopologicaldef}
Let $B$ be a topological space and $p \in \mathbb{Z}_+$. A
\textbf{singular $C_p$-space} in $B$ is a continuous map $\sigma:P 
\rightarrow B$ where $P\in C_p$, and the \textbf{singular $C_p$-chain group} 
$S_p(B)$ is the free abelian group generated by the singular $C_p$-spaces. 
Define $S_p(B) = \{0\}$ if $S_p = \{0\}$ or $B =\emptyset$.
Elements of $S_p(B)$ are called \textbf{singular topological chains}
of \textbf{degree} $p$. 
\end{definition}

\noindent
For $p \geq 1$ there is a boundary operator $\partial_p:S_p(B) \rightarrow
S_{p-1}(B)$ induced from the boundary operator $\partial_p:S_p \rightarrow S_{p-1}$.
If $\sigma:P \rightarrow B$ is a singular $C_p$-space in $B$, then $\partial_p(\sigma)$
is given by the formula
$$
\partial_p(\sigma) = \sum_k n_k \sigma|_{P_k}
$$ 
where 
$$
\partial_p(P) = \sum_k n_k P_k.
$$
The pair $(S_\ast(B),\partial_\ast)$ is called a \textbf{chain complex of
singular topological chains} in $B$.

\smallskip\noindent
Note: The preceeding definitions are quite general.  To construct
the Morse-Bott-Smale chain complex in Section \ref{MBcomplex} we only need $C_p$ 
to include the $p$-dimensional faces of an $N$-cube, the compactified moduli
spaces of gradient flow lines of dimension $p$, and the components of
their fibered products of dimension $p$. In this section we show that these
spaces carry the structure of abstract topological chains.

\subsection{Singular $N$-cube chains}\label{Ncubes}
In this subsection we show how the above abstract definitions apply to the
faces of an $N$-cube to produce a chain complex comparable to the singular cubical chain 
complex of \cite{MasABa}.

\smallskip
Pick some large positive integer $N$ and let $I^N = \{(x_1,\ldots ,x_N) \in
\mathbb{R}^N |\ 0\leq x_j \leq 1, \ j=1,\ldots ,N \}$ denote the unit
$N$-cube.  For every $0 \leq p \leq N$ let $C_p$ be the set consisting of
the faces of $I^N$ of dimension $p$, i.e. subsets of $I^N$ where $p$ of the
coordinates are free and the rest of the coordinates are fixed to be either
$0$ or $1$. For every $0 \leq p \leq N$ let $S_p$ be the free abelian group
generated by the elements of $C_p$.

For $P \in C_p$ we define
$$
\partial_p(P) = \sum_{j=1}^p (-1)^j \left[ P|_{x_j = 1} - P|_{x_j = 0} \right]
\in S_{p-1}
$$
where $x_j$ denotes the $j^{\text{th}}$ free coordinate of $P$. The boundary operator
$\partial_p$ extends linearly to a homomorphism $\partial_p:S_p \rightarrow S_{p-1}$.
To see that $\partial_{p-1} \circ \partial_p = 0$ note that in the sum for 
$\partial_{p-1}(\partial_p( P ))$,
$$
\sum_{i=1}^{p-1} (-1)^i  \left[ \left.
\sum_{j=1}^p (-1)^j \left[ P|_{x_j = 1} - P|_{x_j = 0} \right] \right|_{\tilde{x}_i = 1} -\ \ \ 
\left. \sum_{j=1}^p (-1)^j \left[ P|_{x_j = 1} - P|_{x_j = 0} \right] \right|_{\tilde{x}_i = 0} \right]
$$
the terms cancel in pairs and hence $\partial_{p-1}(\partial_p(P)) = 0 \in S_{p-2}$.
A continuous map $\sigma_P:P \rightarrow B$ into a topological space $B$ is
a singular $C_p$-space in $B$. The boundary operator applied to $\sigma_P$
is
$$
\partial_p(\sigma_P) = \sum_{j=1}^p (-1)^j \left[\sigma_P|_{x_j = 1} -
\sigma_P|_{x_j = 0}\right] \in S_{p-1}(B)
$$
where $\sigma_P|_{x_j = 0}$ denotes the restriction $\sigma_P: P|_{x_j=0}
\rightarrow B$ and $\sigma_P|_{x_j = 1}$ denotes the restriction
$\sigma_P:P|_{x_j = 1} \rightarrow B$. We will show that the $p^\text{th}$ singular 
homology group of $B$ is isomorphic to the homology of a quotient of the chain complex 
$(S_\ast(B),\partial_\ast)$ by the degenerate singular chains for $p < N$.

\begin{definition}[Degeneracy Relations for Singular $N$-Cube Chains]\label{cubedegenerate}
Let $\sigma_P$ and $\sigma_Q$ be singular $C_p$-spaces in $B$ and let
$\partial_p(Q) = \sum_j n_j Q_j \in S_{p-1}$.  For any map $\alpha:P 
\rightarrow Q$, let $\partial_p(\sigma_Q) \circ \alpha$ denote the 
formal sum $\sum_j n_j (\sigma_Q\circ\alpha) |_{\alpha^{-1} (Q_j)}$. 
Define the subgroup $D_p(B)\subseteq S_p(B)$ of \textbf{degenerate singular
$N$-cube chains} to be the subgroup generated by the following elements. 
\begin{enumerate}
\item If $\alpha$ is an orientation preserving homeomorphism such that
      $\sigma_Q \circ \alpha = \sigma_P$ and $\partial_p(\sigma_Q) \circ 
      \alpha = \partial_p(\sigma_P)$, then $\sigma_P - \sigma_Q \in D_p(B)$.
\item If $\sigma_P$ does not depend on some free coordinate of 
      $P$, then $\sigma_P \in D_p(B)$.
\end{enumerate}
\end{definition}

\noindent
Note: We orient the faces of $I^N$ of dimension $p$ by requiring that the map
which identifies the free coordinates $(x_1,\ldots ,x_p)$ of a face with
the standard coordinates of $\mathbb{R}^p$ be orientation preserving.
Also, a map from a vertex to a vertex is orientation preserving.

\begin{theorem}[Singular $N$-Cube Chain Theorem]\label{cubehomology}
The boundary operator for singular $N$-cube chains 
$\partial_p:S_p(B) \rightarrow S_{p-1}(B)$ descends to a homomorphism 
$$
\partial_p:S_p(B)/D_p(B) \rightarrow S_{p-1}(B)/D_{p-1}(B),
$$
and
$$
H_p(S_\ast(B)/D_\ast(B), \partial_\ast) \approx H_p(B;\mathbb{Z})
$$ 
for all $p < N$.
\end{theorem}

\smallskip\noindent
Proof:
Let $\sigma_P \in S_p(B)$ be a singular $C_p$-space and assume that there
exists a singular $C_p$-space $\sigma_Q \in S_p(B)$ and an orientation
preserving homeomorphism $\alpha:P \rightarrow Q$ such that 
$\sigma_Q \circ \alpha = \sigma_P$ and $\partial_p(\sigma_Q) \circ 
\alpha = \partial_p(\sigma_P)$. Since $\partial_{p-1}(\partial_p 
(\sigma_Q)) = \partial_{p-1}(\partial_p (\sigma_P)) = 0$ we have 
$\partial_{p-1}(\partial_p(\sigma_Q)) \circ \alpha = 
\partial_{p-1}(\partial_p (\sigma_P))$. Thus, $\partial_p(\sigma_P - 
\sigma_Q) = \partial_p(\sigma_P) - \partial_p(\sigma_Q) \in D_{p-1}(B)$.

Now assume that $\sigma_P\in S_p(B)$ is a singular $C_p$-space that does
not depend on some free coordinate $x_j$. Then in the sum for
$\partial_p(\sigma_P)$ the term $\sigma_P|_{x_j = 1} - \sigma_P|_{x_j = 0}
\in D_{p-1}(B)$ because $\sigma_P|_{x_j = 0}$ and $\sigma_P|_{x_j = 1}$
satisfy the first condition for degeneracy, and the rest of the terms
are all independent of the coordinate corresponding to $x_j$ and thus
satisfy the second condition for degeneracy. Therefore
$\partial_p(D_p(B)) \subseteq D_{p-1}(B)$, and the boundary operator
$\partial_p$ descends to a homomorphism of the quotient groups.

We now recall the definition of singular homology using singular 
cubes in \cite{MasABa}. Let $I^p$ be the standard unit $p$-cube for 
$p>0$ and let $I^0$ be the origin in $\mathbb{R}^N$. We will view 
$I^p\in C_p$ by identifying $I^p$ with $I^p \times (0,\ldots ,0) \subseteq I^N$.
A \textbf{singular $p$-cube} in $B$ is a continuous map $T:I^p \rightarrow B$, 
i.e. a singular $C_p$-space in $B$ with domain $I^p$. Let $Q_p(B)\subseteq S_p(B)$ 
denote the free abelian group generated by the set of all singular 
$p$-cubes in $B$.  A singular $p$-cube is called degenerate if it 
satisfies the second condition listed in Definition \ref{cubedegenerate}.
Note that the first condition is not needed to define singular homology 
using singular $p$-cubes because all singular $p$-cubes are defined on 
the same domain. Let $DQ_p(B)$ denote the subgroup of $Q_p(B)$ 
generated by the degenerate singular $p$-cubes.

If $T:I^p \rightarrow B$ is a singular $p$-cube, then for $i=1,\ldots ,p$
there are singular $p-1$ cubes defined by the formulas
\begin{eqnarray*}
A_iT(x_1,\ldots ,x_{p-1}) = T(x_1,\ldots ,x_{i-1},1,x_i,\ldots ,x_{p-1})\\
B_iT(x_1,\ldots ,x_{p-1}) = T(x_1,\ldots ,x_{i-1},0,x_i,\ldots ,x_{p-1})
\end{eqnarray*}
and a boundary operator defined by
$$
\tilde{\partial}_p(T) = \sum_{i=1}^p (-1)^i [A_iT - B_iT].
$$
The boundary operator satisfies $\tilde{\partial}_p(DQ_p(B)) \subseteq DQ_{p-1}(B)$,
and the singular homology of $B$ is the homology of the chain complex
$(Q_\ast(B)/DQ_\ast(B),\tilde{\partial}_\ast)$.  For more details see
Chapter VII of \cite{MasABa}.

Now assume that $p < N$.   Since $Q_p(B) \subseteq S_p(B)$ 
and $DQ_p(B) \subseteq D_p(B)$ there are induced maps such that
$$
\xymatrix{
Q_p(B)/DQ_p(B) \ar[r]^{\tilde{\partial}_p} \ar[d] &  Q_{p-1}(B)/DQ_{p-1}(B) \ar[d]\\
Q_p(B)/(D_p(B)\cap Q_p(B)) \ar[r]^{\tilde{\partial}_p} \ar[d] & Q_{p-1}(B)/(D_{p-1}(B)\cap Q_p(B)) \ar[d]\\
S_p(B)/D_p(B) \ar[r]^{\partial_p} & S_{p-1}(B)/D_{p-1}(B)
}
$$
commutes for all $p$. We will show that the induced maps in homology are
isomorphisms.  To see that this is true for the bottom square, note that
for any $P \in C_p$ there is an orientation preserving homeomorphism 
$\alpha^{-1}:I^p \rightarrow P$ defined by sending the fixed
coordinates of $I^p \times (0,\ldots ,0) \subseteq I^N$ to the fixed coordinates
of $P$ and the $j^{\text{th}}$ free coordinate of $I^p$ to the $j^{\text{th}}$
free coordinate of $P$ for all $j=1,\ldots ,p$. For any singular $C_p$-space 
$\sigma_P \in S_p(B)$, $\sigma_Q =\sigma_P\circ \alpha^{-1}:I^p \rightarrow B$ 
is a singular $C_p$-space that satisfies $\sigma_Q \circ \alpha = \sigma_P$ 
and $\partial_p(\sigma_Q) \circ \alpha = \partial_p(\sigma_P)$.  Therefore, 
$[\sigma_P] = [\sigma_Q] \in S_p(B)/D_p(B)$ and we see that every equivalence 
class in $S_p(B)/D_p(B)$ has a representative that is defined on the standard 
unit $p$-cube $I^p$.  Therefore, $S_p(B)/D_p(B) \approx Q_p(B)/(D_p(B)
\cap Q_p(B))$ for all $p$, and the homology groups of $(Q_\ast(B)/(D_\ast(B) 
\cap Q_\ast(B)), \tilde{\partial}_\ast)$ are isomorphic to the homology groups of 
$(S_\ast(B)/D_\ast(B),\partial_\ast)$.

To show that the top square in the above diagram induces isomorphisms of the
homology groups we will show that the homology groups of the chain complex
$$
((D_\ast(B)\cap Q_\ast(B))/DQ_\ast(B), \tilde{\partial}_\ast)
$$ 
are zero for all $p$. Thus, the long exact sequence induced from
$$
0 \rightarrow (D_\ast(B)\cap Q_\ast(B))/DQ_\ast(B) \rightarrow
Q_\ast(B)/DQ_\ast(B) \rightarrow Q_\ast(B)/(D_\ast(B) \cap Q_\ast(B)) \rightarrow 0
$$
shows that the homology groups of $(Q_\ast(B)/(D_\ast(B) \cap Q_\ast(B)), 
\tilde{\partial}_\ast)$ are isomorphic to the homology groups of 
$(Q_\ast(B)/DQ_\ast(B), \tilde{\partial}_\ast)$.

Let $G_p = (D_p(B)\cap Q_p(B))/DQ_p(B)$.  Note that $G_0 = \{0\}$ since 
$D_0(B) \cap Q_0(B) = \{0\}$.  Thus, $\tilde{\partial}_1$ is surjective and 
$H_0(G_\ast,\tilde{\partial}_\ast) = \{0\}$. Now assume that $p > 0$ and
let
$$
\sum_k n_k (\sigma_k \circ \alpha_k - \sigma_k) \in D_p(B) \cap Q_p(B)
$$
be a representative for an element of $\text{ker}(\tilde{\partial}_{p}) 
\subseteq G_p$, i.e.
$$
\tilde{\partial_p} \left(\sum_k n_k (\sigma_k \circ \alpha_k - \sigma_k) \right) = 0
$$
and for every $k$ we have $\partial_p(\sigma_k) \circ \alpha_k = 
\partial_p(\sigma_k \circ \alpha_k)$ where $\alpha_k:I^p \rightarrow I^p$
is an orientation preserving homeomorphism. Extend each $\alpha_k$ 
to an orientation preserving homeomorphism $\tilde{\alpha}_k: 
I^{p+1} \rightarrow I^{p+1}$ by the formula $\tilde{\alpha}_k(x_1,\ldots ,x_p,x_{p+1}) = 
(\alpha_k(x_1,\ldots ,x_p),x_{p+1})$, and extend each $\sigma_k:I^p 
\rightarrow B$ to a map $\tilde{\sigma}_k:I^{p+1} \rightarrow B$ 
such that for every $0 \leq x_{p+1} \leq 1$ we have 
$$
\partial_p\left( \sum_k n_k (\tilde{\sigma}_k(-,x_{p+1})) \circ 
\tilde{\alpha}_k - \tilde{\sigma}_k(-,x_{p+1}) \right) = 0,
$$
$\partial_p(\tilde{\sigma}_k(-,x_{p+1})) \circ \tilde{\alpha}_k = 
\partial_p(\tilde{\sigma}_k \circ \tilde{\alpha}_k(-,x_{p+1}))$, 
and each $\tilde{\sigma}_k$ is constant when $x_{p+1} = 1$. Then 
$$
\sum_k n_k (\tilde{\sigma}_k \circ \tilde{\alpha}_k - \tilde{\sigma}_k)
\in D_{p+1}(B) \cap Q_{p+1}(B)
$$
satisfies
$$
\tilde{\partial}_{p+1}\left[ \sum_k n_k(\tilde{\sigma}_k \circ \tilde{\alpha}_k
- \tilde{\sigma}_k) \right] = (-1)^{p+1}\left[\sum_k n_k 
(\sigma_k \circ \alpha_k - \sigma_k) \right]
$$
in $(D_p(B)\cap Q_p(B))/DQ_p(B)$. Therefore, $\tilde{\partial}_{p+1}$ 
maps onto $\text{ker}(\tilde{\partial}_p)$ and $(G_\ast,\tilde{\partial}_\ast)$
is acyclic.

\begin{flushright}
$\Box$
\end{flushright}


\subsection{The fibered product of singular topological chains}\label{fiberedsub}
In this subsection we resume the discussion of our abstract definitions. We define 
the fibered product of singular topological chains and show that the fibered product 
of singular topological chains is an abstract topological chain.

\smallskip

Let $\sigma_i:P_i \rightarrow B$ for $i=1,2$ be two continuous maps
into a topological space $B$.  Recall that the fibered product of $\sigma_1$ 
and  $\sigma_2$ is defined to be
$$
P_1 \times_B P_2 = (\sigma_1 \times \sigma_2)^{-1}(\Delta)
$$
where $\Delta$ is the diagonal in $B \times B$, i.e.
$$
P_1 \times_B P_2 = \{(x_1,x_2) \in P_1 \times P_2 |\ \sigma_1(x_1) =
\sigma_2(x_2)\}.
$$

\begin{lemma}\label{fiberprod}
Suppose that $\sigma_1:P_1\rightarrow B$ and $\sigma_2:P_2\rightarrow B$ are smooth
maps where $P_1$, $P_2$, and $B$ are smooth manifolds of dimension $p_1$, $p_2$, and $b$
respectively. If $\sigma_1$ is transverse to $\sigma_2$, then the fibered product 
$P_1 \times_{B} P_2$
is a smooth manifold of dimension $p_1+p_2-b$.
\end{lemma}

\smallskip\noindent
Proof: 
This follows from the fact that 
$\sigma_1 \pitchfork \sigma_2$ if and only if $(\sigma_1 \times \sigma_2)
\pitchfork \Delta$.

\begin{flushright}
$\Box$
\end{flushright}

Given a collection of topological spaces $\{C_p\}_{p \geq 0}$ we will say that 
an element $P \in C_p$ has \textbf{degree} $p$. If $B$ is a topological 
manifold of dimension $b$, $P_1 \in C_{p_1}$, $P_2 \in C_{p_2}$, and 
$\sigma_i:P_i \rightarrow B$ for $i=1,2$ are continuous maps, 
then we can associate the \textbf{degree} $p_1 + p_2 - b$ to the fibered product
$P_1 \times_B P_2$ using the same formula as above. If we assume that the
collection of spaces $\{C_p\}_{p \geq 0}$ is closed under the fibered product 
construction with respect to some collection of maps, i.e. if
$$
P_1 \times_B P_2 \in C_{p_1 + p_2 - b},
$$
then the fibered product construction extends linearly to the collection of free 
abelian groups $\{S_p\}$.

\begin{definition}\label{fiberedboundary}
Suppose that $\{C_p\}_{p \geq 0}$ is a collection of topological spaces
that is closed under the fibered product construction with respect
to some collection of maps, and assume that $(S_\ast, \partial_\ast)$ is
a chain complex of abstract topological chains based on
the collection $\{C_p\}_{p \geq 0}$.
If $\sigma_i = \sum_k n_{i,k} \sigma_{i,k} \in S_{p_i}(B)$ for $i=1,2$ 
where $\sigma_{i,k}:P_{i,k} \rightarrow B$ is a singular $C_{p_i}$-space 
for all $k$, then the \textbf{fibered product} of $\sigma_1$ and 
$\sigma_2$ over $B$ is defined to be
$$
P_1 \times_B P_2 = \sum_{k,j} n_{1,k} n_{2,j}\ P_{1,k} \times_B P_{2,j}
$$
where $P_1 = \sum_k n_{1,k}P_{1,k} \in S_{p_1}$ and 
$P_2 = \sum_j n_{2,j}P_{2,j} \in S_{p_2}$. The \textbf{boundary operator} 
applied to the fibered product is defined to be
$$
\partial (P_1 \times_B P_2) = \partial P_1 \times_B P_2
+ (-1)^{p_1+b} P_1 \times_B \partial P_2.
$$
If $\sigma_i = 0$ for either $i=1$ or $2$, then we define $P_1 \times_B P_2 = 0$.
\end{definition}

\noindent
Note: We omit the subscript $p$ on the boundary operator 
$\partial_p$ in order to simplify the notation.  Also, given the
data of a triple
$$
\xymatrix{
P_1 \ar[r]^{\sigma_{11}} & B_1 & P_2 \ar[l]_{\sigma_{12}} \ar[r]^{\sigma_{22}} & B_2 & \ar[l]_{\sigma_{23}} P_3\\
}
$$
we can form the iterated fibered product $(P_1 \times_{B_1} P_2) \times_{B_2} P_3$
using $\sigma_{23}$ and the map $\sigma_{22} \circ \pi_2:P_1 \times_{B_1} P_2
\rightarrow B_2$, where $\pi_2:P_1 \times_{B_1} P_2 \rightarrow P_2$
denotes projection to the second component.  That is, we have the
following diagram.
$$
\xymatrix{
(P_1 \times_{B_1} P_2) \times_{B_2} P_3 \ar@{-->}[d]^-{\pi_1} \ar@{-->}[rr]^-{\pi_3} &  & P_3 \ar[d]^{\sigma_{23}}\\
P_1 \times_{B_1} P_2 \ar@{-->}[d]^-{\pi_1} \ar@{-->}[r]^-{\pi_2} & P_2 \ar[d]^{\sigma_{12}} \ar[r]^{\sigma_{22}} & B_2 \\
P_1 \ar[r]^{\sigma_{11}} & B_1 & \\
}
$$
Similarly, we can form the iterated fibered product $P_1 \times_{B_1} (P_2 \times_{B_2} P_3)$
using $\sigma_{11}$ and the map $\sigma_{12} \circ \pi_1:P_1 \times_{B_1} P_2
\rightarrow B_1$, where $\pi_1:P_2 \times_{B_2} P_3 \rightarrow P_2$
denotes projection to the first component.  This corresponds to the
following diagram.
$$
\xymatrix{
P_1 \times_{B_1} (P_2 \times_{B_2} P_3) \ar@{-->}[dd]^-{\pi_1} \ar@{-->}[r]^-{\pi_2} & P_2 \times_{B_2} P_3 \ar@{-->}[d]^-{\pi_1} 
   \ar@{-->}[r]^-{\pi_2} & P_3 \ar[d]^{\sigma_{23}}\\
 & P_2 \ar[d]^{\sigma_{12}} \ar[r]^{\sigma_{22}} & B_2 \\
P_1 \ar[r]^{\sigma_{11}} & B_1 & \\
}
$$
It is easy to see that $(P_1 \times_{B_1} P_2) \times_{B_2} P_3$
and $P_1 \times_{B_1} (P_2 \times_{B_2} P_3)$ are the same as topological
spaces.  The following lemma shows that they are also the same as 
abstract topological chains.

\begin{lemma}\label{fiberedcomplex}
The fibered product of two singular topological chains is an abstract
topological chain, i.e. the boundary operator on fibered products is of
degree -1 and satisfies $\partial \circ \partial = 0$. Moreover, 
the boundary operator on fibered products is associative, i.e.  
$$
\partial((P_1 \times_{B_1} P_2) \times_{B_2} P_3) = \partial(P_1 \times_{B_1} (P_2 \times_{B_2} P_3)).
$$
\end{lemma}

\smallskip\noindent
Proof: Since $\partial$ is a boundary operator on $P_1$ and $P_2$, the degree of
$\partial P_1$ is $p_1 -1$ and the degree of $\partial P_2$ is $p_2 - 1$. 
Hence both $\partial P_1 \times_B P_2$ and $P_1 \times_B \partial P_2$
have degree $p_1 + p_2 - b - 1$.

To see that $\partial^2 (P_1 \times_B P_2) = 0$ we compute as follows.

\begin{eqnarray*}
\partial(\partial( P_1 \times_B P_2)) & = & \partial( \partial P_1 \times_B P_2 +
(-1)^{p_1+b} P_1 \times_B \partial P_2 )\\
& = & \partial^2 P_1 \times_B P_2 + (-1)^{p_1 - 1 + b} \partial P_1 \times_B
      \partial P_2 +\\
&   & \quad (-1)^{p_1 + b} (\partial P_1 \times_B \partial P_2 + (-1)^{p_1 + b}
      P_1 \times_B \partial^2 P_2) \\
& = & 0.
\end{eqnarray*}

To prove associativity we compute as follows.
\smallskip

$\partial(P_1 \times_{B_1} (P_2 \times_{B_2} P_3))$
$$
\begin{array}{ll}
= & \partial P_1 \times_{B_1} (P_2 \times_{B_2} P_3) + (-1)^{p_1+b_1} P_1
    \times_{B_1}\partial (P_2 \times_{B_2} P_3)\\
= & \partial P_1 \times_{B_1} (P_2 \times_{B_2} P_3) +\\
  & \quad (-1)^{p_1 + b_1} (P_1 \times_{B_1} (\partial P_2 \times_{B_2} P_3 +
    (-1)^{p_2 + b_2}P_2 \times_{B_2} \partial P_3 ))\\
= & \partial P_1 \times_{B_1} P_2 \times_{B_2} P_3 + (-1)^{p_1 + b_1} P_1
    \times_{B_1}\partial P_2 \times_{B_2} P_3 +\\
  & \quad (-1)^{p_1+p_2+b_1+b_2} P_1 \times_{B_1} P_2 \times_{B_2} \times 
    \partial P_3 
\end{array}
$$

\smallskip

$\partial((P_1 \times_{B_1} P_2) \times_{B_2} P_3)$
$$
\begin{array}{ll}
= & \partial (P_1 \times_{B_1} P_2) \times_{B_2} P_3 + (-1)^{\mbox{deg}(P_1
    \times_{B_1} P_2) +b_2}(P_1 \times_{B_1} P_2) \times_{B_2} \partial P_3\\
= & (\partial P_1 \times_{B_1} P_2 + (-1)^{p_1+b_1}P_1 \times_{B_1} \partial P_2) 
    \times_{B_2} P_3 + \\
  & \quad (-1)^{p_1+p_2-b_1+b_2} P_1 \times_{B_1} P_2 \times_{B_2} \partial P_3\\
= & \partial P_1 \times_{B_1} P_2 \times_{B_2} P_3 + (-1)^{p_1 + b_1} P_1 \times_{B_1}
    \partial P_2 \times_{B_2} P_3 +\\
  & \quad (-1)^{p_1+p_2-b_1+b_2} P_1 \times_{B_1} P_2 \times_{B_2} \times \partial P_3 
\end{array}
$$

\begin{flushright}
$\Box$
\end{flushright}


\subsection{Compactified moduli spaces as singular topological chains}\label{modulisingular}
In this subsection we define a degree and boundary operator for compactified 
moduli spaces of gradient flow lines. We will show in Subsection \ref{manicorners}
that these compactified moduli spaces are smooth manifolds with corners and
the degree defined in this subsection coincides with the dimension of these manifolds.
The abstract topological chain structure defined in this subsection together 
with those defined for the faces of an $N$-cube and for fibered products are 
used in Section \ref{MBcomplex} to define the homomorphism $\partial_0$ in 
the Morse-Bott-Smale chain complex.

\medskip
Let $f:M \rightarrow \mathbb{R}$ be a Morse-Bott function on a Riemannian
manifold $M$ and let $\varphi_t:M \rightarrow M$ denote the flow associated
to $-\nabla f$.  For any two critical submanifolds $B$ and $B'$ the flow
$\varphi_t$ induces an $\mathbb{R}$-action on $W^u(B) \cap W^s(B')$.
Let
$$
\mathcal{M}(B,B') = (W^u(B) \cap W^s(B'))/\mathbb{R}
$$
be the quotient space of gradient flow lines from $B$ to $B'$.  For proofs of
the following two fundamental theorems we refer the reader to Appendix \S A.3
of \cite{AusMor} (see \cite{BouCom}, \cite{HofAGeI}, and \cite{HofAGeII} for a different 
approach to Theorem \ref{compactification}).

\begin{theorem}[Gluing]\label{gluing}
Suppose that $B$, $B'$, and $B''$ are critical submanifolds such that
$W^u(B) \pitchfork W^s(B')$ and $W^u(B') \pitchfork W^s(B'')$. In addition,
assume that $W^u(x) \pitchfork W^s(B'')$ for all $ x \in B'$. Then
for some $\epsilon>0$, there is an injective local diffeomorphism
$$
G:\mathcal{M}(B,B') \times_{B'}
\mathcal{M}(B',B'') \times (0,\epsilon) \rightarrow \mathcal{M}(B,B'')
$$
onto an end of $\mathcal{M}(B,B'')$.
\end{theorem}

\begin{theorem}[Compactification]\label{compactification}
Assume that $f:M \rightarrow \mathbb{R}$ satisfies the Morse-Bott-Smale
transversality condition. For any two distinct connected critical submanifolds $B$
and $B'$ the moduli space $\mathcal{M}(B,B')$ has a compactification 
$\overline{\mathcal{M}}(B,B')$, consisting of all the piecewise gradient
flow lines from $B$ to $B'$, which is either empty or a compact smooth manifold with 
corners of dimension $\lambda_B - \lambda_{B'}+b-1$. Moreover, the beginning and endpoint
maps extend to smooth maps
$$
\begin{array}{l}
\partial_- : \overline{\mathcal{M}}(B,B') \rightarrow B \\
\partial_+ : \overline{\mathcal{M}}(B,B') \rightarrow B',
\end{array}
$$
where the beginning point map $\partial_-$ has the structure of a locally trivial fiber bundle.
\end{theorem}

The compactified moduli spaces of gradient flow lines for a Morse-Bott-Smale
function can be described inductively as follows. For any two critical 
submanifolds $B$ and $B'$, define $B \succ B'$ if $B \neq B'$ and there 
exists a gradient flow line from $B$ to $B'$. Note that this relationship 
is transitive, i.e. if $B \succ B'$ and $B' \succ B''$, then $B \succ B''$ 
by Theorem \ref{gluing}. If there is no critical submanifold $B'$ such 
that $B \succ B' \succ B''$, then $\overline{\mathcal{M}}(B,B'') = 
\mathcal{M}(B,B'')$.  In general,
$$
\overline{\mathcal{M}}(B,B'') = \mathcal{M}(B,B'') \cup \bigcup_{B \succ B'
\succ B''} \overline{\mathcal{M}}(B,B') \times_{B'}
\overline{\mathcal{M}}(B',B'')
$$
where the union is taken over all critical submanifolds $B'$ between
$B$ and $B''$.  Note that this description is inductive, i.e. both
$\overline{\mathcal{M}}(B,B')$ and $\overline{\mathcal{M}}(B',B'')$ may contain
fibered products of compactified moduli spaces of smaller dimension. Hence
$\overline{\mathcal{M}}(B,B'')$ consists of all piecewise gradient flow lines from
$B$ to $B''$. The topology of $\overline{\mathcal{M}}(B,B')$ is determined 
by considering its elements as one dimensional subsets of $M$.
\begin{definition}
The topology on $\overline{\mathcal{M}}(B,B')$ is the topology induced
from the Hausdorff metric $d_H$ on subsets of $M$.  That is,
if $l_1,\ l_2 \in \overline{\mathcal{M}}(B,B')$, then $l_1$ and $l_2$
determine subsets $l_1,l_2 \subseteq M$, which are images of continuous 
injective paths from $B$ to $B'$, and 
$$
d_H(l_1,l_2) = \sup_{x_1\in l_1} \inf_{x_2\in l_2} d(x_1,x_2) +
\sup_{x_2\in l_2} \inf_{x_1\in l_1} d(x_1,x_2)
$$
where $d$ is the metric on $M$.
\end{definition}
\smallskip\noindent
Note that this definition is compatible with the topology of
fibered products. That is, if $B\succ B' \succ B''$, then 
$\overline{\mathcal{M}}(B,B')\times_{B'} \overline{\mathcal{M}}(B',B'')$
inherits a topology as a subspace of $\overline{\mathcal{M}}(B,B') \times 
\overline{\mathcal{M}}(B',B'')$, and this topology coincides with the 
Hausdorff topology on subsets of $M$ \cite{HurFlo}.

The next definition shows how to use the above description of the compactified
moduli spaces to define a degree and boundary operator satisfying the axioms
for abstract topological chains. For this definition we only need to assume that
the Morse-Bott function $f:M \rightarrow \mathbb{R}$ is weakly self-indexing
(see Lemma \ref{weakly}). To simplify the notation we will assume that for
each $i=0,\ldots ,m$ the components of $B_i$, the set of critical points
of index $i$, are of the same dimension. In general one needs to group the 
components by their dimension and then define the degree and boundary operator 
on each group.

\begin{definition}\label{moduliboundary}
Let $B_i$ be the set of critical points of index $i$. For any $j=1,\ldots ,i$
the degree of $\overline{\mathcal{M}}(B_i,B_{i-j})$ is defined to be $j+b_i -1$
and the boundary operator is defined to be
$$
\partial \overline{\mathcal{M}}(B_i,B_{i-j}) = (-1)^{i+b_i} \sum_{i-j<n<i}
\overline{\mathcal{M}}(B_i,B_n) \times_{B_n} \overline{\mathcal{M}}(B_n,B_{i-j})
$$
where $b_i = \mbox{dim }B_i$ and the fibered product is taken over the beginning
and endpoint maps $\partial_-$ and $\partial_+$. If $B_n = \emptyset$, then
$\overline{\mathcal{M}}(B_i,B_n) = \overline{\mathcal{M}}(B_n,B_{i-j}) = 0$.
The boundary operator extends to fibered products of compactified moduli spaces 
via Definition
\ref{fiberedboundary}.
\end{definition}

\noindent
Note: For all $p\geq 0$ let $C_p$ be the set consisting of
the connected components of degree $p$ of fibered products of the form
$$
\overline{\mathcal{M}}(B_{i_1},B_{i_2}) \times_{B_{i_2}} 
\overline{\mathcal{M}}(B_{i_2},B_{i_3}) \times_{B_{i_3}} \cdots \times_{B_{i_{n-1}}}
\overline{\mathcal{M}}(B_{i_{n-1}},B_{i_n})
$$
where $m \geq i_1 > i_2 > \cdots > i_n \geq 0$ and the fibered products are 
taken with respect to the beginning and endpoint maps $\partial_-$ and
$\partial_+$. Let $S_p$ be the free abelian group generated by the
elements of $C_p$, i.e. $S_p = \mathbb{Z}[C_p]$.  The preceding
definition describes $\overline{\mathcal{M}}(B_i,B_{i-j})$ and
$\partial \overline{\mathcal{M}}(B_i,B_{i-j})$ as abstract topological 
chains in $S_{j+b_i-1}$ and $S_{j+b_i-2}$ respectively.  If $\overline{\mathcal{M}}(B_i,B_{i-j})$
has more than one connected component, then $\overline{\mathcal{M}}(B_i,B_{i-j})
\in S_{j+b_i-1}$ is defined to be the sum of these components.
Similarly, $\overline{\mathcal{M}}(B_i,B_n) \times_{B_n} \overline{\mathcal{M}}(B_n,B_{i-j})
\in S_{j+b_i-2}$ is defined to be the sum of its connected
components.

\noindent
Notation: We will use the following notation for the fibered product over the
beginning and endpoint maps $\partial_-$ and $\partial_+$.  For any $i-j < n < i$ we
define
$$
\overline{\mathcal{M}}(B_i,B_n,B_{i-j}) =
\overline{\mathcal{M}}(B_i,B_n) \times_{B_n} \overline{\mathcal{M}}(B_n,B_{i-j}).
$$
Similarly for $i-j<n<s<i$, $\overline{\mathcal{M}}(B_i,B_s,B_n,B_{i-j})$
will denote the triple fibered product over $B_s$ and $B_n$. Iterated fibered
products are well defined as abstract topological chains since the boundary
operator on fibered products is associative by Lemma \ref{fiberedcomplex}.

\begin{lemma}\label{modulicomplex}
The degree and boundary operator for $\overline{\mathcal{M}}(B_i,B_{i-j})$
satisfy the axioms for abstract topological chains, i.e. the boundary operator
on compactified moduli spaces of gradient flow lines is of degree $-1$
and it satisfies $\partial \circ \partial = 0$.
\end{lemma}

\smallskip\noindent
Proof: For any $n$ with $i-j<n<i$ the degree of $\overline{\mathcal{M}}(B_i,B_n)$
is $i-n+b_i-1$ and the degree of $\overline{\mathcal{M}}(B_n,B_{i-j})$ is
$n-i+j+b_n-1$ by Definition \ref{moduliboundary}.  Thus the degree of
$\overline{\mathcal{M}}(B_i,B_n) \times_{B_n} \overline{\mathcal{M}}(B_n,B_{i-j})$
is $j+b_i-2$ by Definition \ref{fiberedboundary}. We can also apply 
Definition \ref{fiberedboundary} and Definition \ref{moduliboundary}
as follows.  To simplify the notation let $d = \mbox{deg } 
\overline{\mathcal{M}}(B_i,B_n) = i-n + b_i - 1$.  

\medskip\noindent
$
\partial(\overline{\mathcal{M}}(B_i,B_n)\times_{B_n} \overline{\mathcal{M}}(B_n,B_{i-j}))
$
$$
\begin{array}{ll}
= & \displaystyle \partial \overline{\mathcal{M}}(B_i,B_n) \times_{B_n}
    \overline{\mathcal{M}}(B_n,B_{i-j}) + (-1)^{d + b_n} \overline{\mathcal{M}}(B_i,B_n)
    \times_{B_n} \partial \overline{\mathcal{M}}(B_n,B_{i-j}) \\
\\
= & \displaystyle (-1)^{i+ b_i} \sum_{n<s<i} \overline{\mathcal{M}}(B_i,B_s,B_n,B_{i-j})+
    (-1)^{i+b_i-1} \sum_{i-j<t<n} \overline{\mathcal{M}}(B_i,B_n,B_t,B_{i-j})\\
\end{array}
$$

Therefore,
$$
\begin{array}{lll}
\partial^2 \overline{\mathcal{M}}(B_i,B_{i-j})
& = & \displaystyle (-1)^{i+b_i} \left[\sum_{i-j<n<i} \left( (-1)^{i+ b_i}
      \sum_{n<s<i} \overline{\mathcal{M}}(B_i,B_s,B_n,B_{i-j})\ + \right.\right.\\
&   & \displaystyle \hspace{1.3 in} \left.\left.(-1)^{i+b_i-1} \sum_{i-j<t<n}
      \overline{\mathcal{M}}(B_i,B_n,B_t,B_{i-j})\right)\right]\\
& = & \displaystyle (-1)^{i+b_i} \left[(-1)^{i+b_i} \sum_{i-j<n<s<i}  
      \overline{\mathcal{M}}(B_i,B_s,B_n,B_{i-j})\ + \right.\\
&   & \displaystyle \hspace{0.67 in} \left.(-1)^{i+b_i-1} \sum_{i-j<t<n<i}
      \overline{\mathcal{M}}(B_i,B_n,B_t,B_{i-j}) \right]\\
& = & 0
\end{array}
$$

\begin{flushright}
$\Box$
\end{flushright}


\section{The Morse-Bott-Smale chain complex}\label{MBcomplex}

In this section we define the Morse-Bott-Smale chain complex $(C_\ast(f),\partial)$.
Throughout this section we will assume that $f:M \rightarrow \mathbb{R}$ is 
a Morse-Bott-Smale function on a compact oriented smooth Riemannian manifold
$(M,g)$ of dimension $m$. Furthermore, we will assume that all the critical 
submanifolds $B$ and their negative normal bundles $\nu^-_\ast(B)$
are oriented. If $f$ is constant, then the Morse-Bott-Smale chain complex 
reduces to the chain complex of smooth singular $N$-cube chains (Example
\ref{constantfunction}). If $f$ is a Morse-Smale function, 
then the Morse-Bott-Smale chain complex reduces to the Morse-Smale-Witten 
chain complex (Example \ref{MSfunction}).


\subsection{The Morse-Bott degree and the chain complex $(\tilde{C}_\ast(f),\partial)$}\label{MBchaindef}

Recall that $B_i$ is the set of critical points of index $i$ for
$i=0,\ldots ,m$ and let $N > \text{dim }M$. For any $p\geq 0$ let
$C_p$ be the set consisting of the faces of $I^N$ of dimension $p$ and the 
connected components of degree $p$ of fibered products of the form
$$
Q \times_{B_{i_1}} \overline{\mathcal{M}}(B_{i_1},B_{i_2}) \times_{B_{i_2}} 
\overline{\mathcal{M}}(B_{i_2},B_{i_3}) \times_{B_{i_3}} \cdots \times_{B_{i_{n-1}}}
\overline{\mathcal{M}}(B_{i_{n-1}},B_{i_n})
$$
where $m \geq i_1 > i_2 > \cdots > i_n \geq 0$, $Q$ is a face of $I^N$ of
dimension $q\leq p$, $\sigma:Q \rightarrow B_{i_1}$ is smooth, and the
fibered products are taken with respect to $\sigma$ and the beginning and
endpoint maps $\partial_-$ and $\partial_+$.  The following lemma is 
proved at the end of this section.

\begin{lemma}\label{basiscorners}
The elements of $C_p$ are compact oriented smooth manifolds with corners.
\end{lemma}

\noindent
Let $S_p$ be the free abelian group generated by the elements of $C_p$, and let 
$S_p^\infty(B_i)$ denote the subgroup of the singular $C_p$-chain group $S_p(B_i)$
generated by those maps $\sigma:P \rightarrow B_i$ that satisfy the 
following two conditions: 
\begin{enumerate}
\item The map $\sigma$ is smooth.
\item If $P\in C_p$ is a connected component of a fibered product, 
      then $\sigma = \partial_+ \circ \pi$, where $\pi$ denotes projection 
      onto the last component of the fibered product.
\end{enumerate}

\begin{definition}\label{MBdegree}
Define the \textbf{Morse-Bott degree} of the singular topological chains
in $S_p^\infty(B_i)$ to be $p+i$. For any $k=0,\ldots ,m$ the group 
of smooth singular topological chains of Morse-Bott degree $k$ is defined to be
$$
\tilde{C}_k(f) = \bigoplus_{i=0}^m S_{k-i}^\infty(B_i).
$$
\end{definition}

If $\sigma:P \rightarrow B_i$ is a singular $C_p$-space in $S_p^\infty(B_i)$,
then for any $j = 1,\ldots ,i$ composing the projection map $\pi_2$ onto the second 
component of $P \times_{B_i} \overline{\mathcal{M}}(B_i,B_{i-j})$ with the 
endpoint map $\partial_+:\overline{\mathcal{M}}(B_i,B_{i-j}) \rightarrow B_{i-j}$ 
gives a map
$$
P \times_{B_i} \overline{\mathcal{M}}(B_i,B_{i-j}) 
\stackrel{\pi_2}{\longrightarrow} \overline{\mathcal{M}}(B_i,B_{i-j})
\stackrel{\partial_+}{\longrightarrow} B_{i-j}.
$$
The next lemma shows that restricting this map to the connected components of 
the fibered product $P \times_{B_i}\overline{\mathcal{M}} (B_i,B_{i-j})$ and 
adding these restrictions (with the sign determined by the orientation 
when the dimension of a component is zero) defines an element 
$\partial_j(\sigma) \in S_{p+j-1}^\infty(B_{i-j})$.

\smallskip
\begin{lemma}\label{partialj}
If $\sigma:P \rightarrow B_i$ is a singular $C_p$-space in $S_p^\infty(B_i)$,
then for any $j=1,\ldots ,i$ adding the components of $P \times_{B_i} 
\overline{\mathcal{M}}(B_i,B_{i-j})$ (with sign when the dimension of a component
is zero) yields an abstract topological chain of degree $p+j-1$. That is, we can identify
$$
P \times_{B_i} \overline{\mathcal{M}}(B_i,B_{i-j}) \in S_{p+j-1}.
$$
Thus, for all $j=1,\ldots ,i$ there is an induced homomorphism
$$
\partial_j:S_p^\infty(B_i) \rightarrow S_{p+j-1}^\infty(B_{i-j})
$$
which decreases the Morse-Bott degree by $1$.
\end{lemma}

\smallskip\noindent
Proof:  The set $B_i$ is a union of submanifolds of $M$, possibly of 
different dimensions. Since $\sigma:P \rightarrow B_i$ is continuous and $P$
is connected, the image $\sigma(P)$ must lie in a connected component 
$B \subseteq B_i$.  By Lemma \ref{modulicomplex}, $\overline{\mathcal{M}}(B,B_{i-j})$
is an abstract topological chain of degree $j+b-1$, and by the discussion
after Lemma \ref{fiberprod}, $P \times_{B} \overline{\mathcal{M}}(B,B_{i-j})$ 
has degree $p+j-1$.  This degree is independent of the dimension of the
connected component $B$.

To see that $P \times_{B_i} \overline{\mathcal{M}}(B_i,B_{i-j})$ has finitely
many components, first note that $P \in C_p$ is a compact smooth
manifold with corners by Lemma \ref{basiscorners}, and if $B'$ is a 
connected component of $B_{i-j}$, then $\overline{\mathcal{M}}(B,B')$ 
is a compact smooth manifold with corners by Theorem \ref{compactification}.
We will show in Lemma \ref{fiberedmanifold} that 
$$
P \times_B \overline{\mathcal{M}}(B,B')
$$
is a compact smooth manifold with corners, and hence it has finitely
many components. By assumption (see Definition \ref{MorseBottDefinition}),
the set $B_{i-j}$ has finitely many components. Thus,
$$
P \times_B     \overline{\mathcal{M}}(B,B_{i-j}) = 
P \times_{B_i} \overline{\mathcal{M}}(B_i,B_{i-j})
$$
has finitely many components.

\begin{flushright}
$\Box$
\end{flushright}

\begin{definition}\label{boundary}
For $k=1,\ldots ,m$ define a homomorphism $\partial:\tilde{C}_k(f) 
\rightarrow \tilde{C}_{k-1}(f)$ as follows. If $\sigma \in S_{p}^\infty(B_i)$
is a singular $S_{p}$-space of $B_i$ where $p=k-i$, then
$$
\partial(\sigma) = \bigoplus_{j=0}^m \partial_j(\sigma)
$$
where $\partial_0$ is $(-1)^{k}$ times the boundary operator on singular
topological chains defined in Section \ref{topologicalchains},
$\partial_j(\sigma) = \partial_+ \circ \pi_2: P \times_{B_i} \overline{\mathcal{M}}(B_i,B_{i-j})
\rightarrow B_{i-j}$ for $j=1,\ldots ,i$, and $\partial_j(\sigma) = 0$ otherwise.
The map $\partial$ extends to a homomorphism
$$
\partial:\bigoplus_{i=0}^m S_{k-i}^\infty(B_i) \longrightarrow
\bigoplus_{i=0}^m S_{k-1-i}^\infty(B_i).
$$
\end{definition}

\noindent
When $m=2$ this homomorphism can be pictured as follows.

$$
\xymatrix{
S_0^\infty(B_2) \ar[r]^{\partial_0} \ar[dr]^{\partial_1} \ar[ddr]|(.33){\partial_2} |!{[d];[dr]}\hole & 0 & & \\
S_1^\infty(B_1)\ar@{}[u]|{\oplus} \ar[r]^(.55){\partial_0} \ar[dr]^{\partial_1} & S_0^\infty(B_1) \ar@{}[u]|{\oplus} \ar[r]^{\partial_0} 
 \ar[dr]^{\partial_1} & 0 & \\
S_2^\infty(B_0) \ar@{}[u]|{\oplus} \ar[r]^{\partial_0} & S_1^\infty(B_0) \ar@{}[u]|{\oplus} \ar[r]^{\partial_0} & S_0^\infty(B_0) \ar@{}[u]|{\oplus}  
 \ar[r]^-{\partial_0} & 0\\
\tilde{C}_2(f) \ar@{}[u]|{\|} \ar[r]^{\partial} & \tilde{C}_1(f) \ar@{}[u]|{\|} \ar[r]^{\partial} & \tilde{C}_0(f) \ar@{}[u]|{\|} \ar[r]^-{\partial} & 0
}
$$

\begin{proposition}\label{complex}
For every $j = 0, \ldots ,m$
$$
\sum_{q=0}^j \partial_q \partial_{j-q} = 0.
$$
\end{proposition}

\smallskip\noindent
Proof:
The case $j=0$ follows from Lemma \ref{fiberedcomplex} and Lemma
\ref{modulicomplex}.  Now, let $\sigma\in S_p^\infty(B_i)$ be a singular
$C_p$-space of $B_i$. Since $\partial_q(\partial_{j-q}(\sigma)) = 0$ if $j >i$
we will assume that $1 \leq j \leq i$. To simplify the notation we give the
following computation in terms of abstract topological chains. This is
sufficient because in $\tilde{C}_\ast(f)$ the only map 
allowed from these fibered products is the endpoint map 
$\partial_+$ composed with the projection map onto the last component. 
We will also assume that for every $i = 0,\ldots ,m$ all the components 
of $B_i$ are of the same dimension $b_i$; in general, the components need 
to be grouped by their dimension and the computation repeated on each group.

When $q=0$ we use Definitions \ref{fiberedboundary}, \ref{moduliboundary},
and \ref{boundary} to compute as follows.
\begin{eqnarray*}
\partial_0(\partial_j(P))& = & \partial_0\left(P \times_{B_i}
\overline{\mathcal{M}}(B_i,B_{i-j})\right)\\
& = & (-1)^{p+i-1} \left(\partial P \times_{B_i} \overline{\mathcal{M}}(B_i,B_{i-j}) +
      (-1)^{p+b_i} P \times_{B_i} \partial \overline{\mathcal{M}}(B_i,B_{i-j})\right)\\
& = & (-1)^{p+i-1} \partial P \times_{B_i} \overline{\mathcal{M}}(B_i,B_{i-j})\ + \\
&   & \quad (-1)^{2p+2b_i+2i-1} \sum_{i-j<n<i} P \times_{B_i} 
      \overline{\mathcal{M}}(B_i,B_n) \times_{B_n} \overline{\mathcal{M}}
      (B_n,B_{i-j})
\end{eqnarray*}

\noindent
If $1 \leq q \leq j-1$, then by Definition \ref{boundary} 
$$
\partial_q(\partial_{j-q}(P)) = P \times_{B_i} \overline{\mathcal{M}}
(B_i,B_{i-j+q})\times_{B_{i-j+q}} \overline{\mathcal{M}}(B_{i-j+q},B_{i-j})
$$
\noindent
and if $q=j$, then
$$
\partial_j(\partial_0(P)) = (-1)^{p+i} \partial P \times_{B_i}
\overline{\mathcal{M}}(B_i,B_{i-j}).
$$
Summing these expressions gives the desired result.

\begin{flushright}
$\Box$
\end{flushright}

\begin{corollary}\label{Cfischain}
The pair $(\tilde{C}_\ast(f),\partial)$ is a chain complex, i.e.
$\partial \circ \partial = 0$.
\end{corollary}


\subsection{Orientations}\label{orientationssub}

In this subsection we describe an explicit set of orientations on the elements of $C_p$.

\smallskip
Recall the assumption that every critical submanifold $B$ and every negative
normal bundle $\nu^-_\ast(B)$ are oriented. For any $p \in B$, the relation
$$
T_pM = T_pB \oplus \nu^+_p(B) \oplus \nu^-_p(B)
$$
determines an orientation on $\nu^+_p(B)$. The stable and
unstable manifolds are oriented by requiring that the injective immersions
$E^+:\nu^+_\ast(B) \rightarrow W^s(B)$ and $E^-:\nu^-_\ast(B) \rightarrow 
W^u(B)$ are orientation preserving (see Theorem \ref{stablemanifold}).
If $N \subseteq M$ is an oriented submanifold, then the normal bundle of 
$N$ is oriented by the relation $T_x(N) \oplus \nu_x(N) = T_x(M)$ for all
$x \in N$. For any two connected critical submanifolds $B$ and $B'$, 
the orientation on $W(B,B') = W^u(B) \pitchfork W^s(B')$ is determined 
by the relation
$$
T_x(M) = T_xW(B,B') \oplus \nu_x(W^s(B')) \oplus \nu_x(W^u(B))
$$
for all $x \in W(B,B')$. Picking a non-critical value $a$ between $f(B')$ 
and $f(B)$ we can identify $\mathcal{M}(B,B') = f^{-1}(a) \cap W(B,B')$. 
An orientation on $\mathcal{M}(B,B')$ is then determined by
$$
T_xW(B,B')= \mbox{span}((-\nabla f)(x))  \oplus T_x\mathcal{M}(B,B')
$$
for all $x \in f^{-1}(a) \cap W(B,B')$. This determines an orientation
on the compact manifold with boundary $\overline{\mathcal{M}}(B,B')$ 
(see for instance Section VI.9 of \cite{BreTop}).

\begin{definition}\label{fiberedorientation}
Suppose that $B$ is an oriented smooth manifold without boundary and 
$P_1$ and $P_2$ are oriented smooth manifolds with corners.  If 
$\sigma_1:P_1 \rightarrow B$ and $\sigma_2:P_2 \rightarrow B$ are 
smooth maps that intersect transversally and stratum transversally, 
then the orientation on the smooth manifold with corners 
$P_1 \times_B P_2$ is defined by the relation
$$
(-1)^{(\mbox{dim }B)(\mbox{dim } P_2)} T_\ast(P_1 \times_B P_2) \oplus
(\sigma_1 \times \sigma_2)^\ast (\nu_\ast(\Delta(B))) = T_\ast( P_1 \times P_2),
$$
where $\nu_\ast(\Delta(B))$ denotes the normal bundle of the diagonal in
$B \times $B.
\end{definition}
\noindent

\noindent
Note: Lemma \ref{fiberprod} extends to the category of smooth manifolds with
corners under some additional assumptions; see the end of this section
for more details. Also, the normal bundle of the diagonal $\Delta(B) \subset B \times B$ 
pulls back to the normal bundle of the fibered product $P_1 \times_{B} P_2$ 
in $P_1 \times P_2$ via $\sigma_1 \times \sigma_2:P_1 \times P_2 \rightarrow
B \times B$.

\begin{lemma}\label{fiberedassociative}
The above orientation on fibered products of transverse intersections of smooth
manifolds with corners is associative, i.e.
$$
(P_1 \times_{B_1} P_2) \times_{B_2} P_3 = P_1 \times_{B_1}(P_2 \times_{B_2} P_3)
$$
as oriented smooth manifolds with corners.
\end{lemma}

\smallskip\noindent
Proof: The proof is a straightforward computation.  The sign
$(-1)^{(\mbox{dim }B)(\mbox{dim } P_2)}$ is needed to prove associativity.

\begin{flushright}
$\Box$
\end{flushright}

We will see in the proof of Lemma \ref{basiscorners} that Definition 
\ref{fiberedorientation} applies to the fibered products in $C_p$, 
and hence by Lemma \ref{fiberedassociative} there is a well-defined orientation 
on the components of 
$$
Q \times_{B_{i_1}} \overline{\mathcal{M}}(B_{i_1},B_{i_2}) \times_{B_{i_2}} 
\overline{\mathcal{M}}(B_{i_2},B_{i_3}) \times_{B_{i_3}} \cdots \times_{B_{i_{n-1}}}
\overline{\mathcal{M}}(B_{i_{n-1}},B_{i_n})
$$
in $C_p$.

\smallskip\noindent
Note: The boundary of $\overline{\mathcal{M}}(B,B')$ consists of fibered
products of compactified moduli spaces of gradient flow lines by 
Theorem \ref{gluing}. These fibered products are oriented by Definition 
\ref{fiberedorientation}. The boundary of $\overline{\mathcal{M}}(B,B')$ 
also inherits an orientation from the orientation on $\overline{\mathcal{M}}(B,B')$ 
(see for instance Lemma VI.9.1 of \cite{BreTop}). These two orientations 
can be compared using Proposition 2.7 of \cite{AusMor} which says that 
the gluing map is orientation reversing.  We will always use the orientation
given by Definition \ref{fiberedorientation}.


\subsection{Degenerate and non-degenerate singular topological chains}\label{degMB}

The Morse-Bott-Smale chain complex $(C_\ast(f),\partial)$ is defined
as the quotient of the chain complex $(\tilde{C}_\ast(f),\partial)$ by the 
degenerate singular topological chains.

\medskip

\begin{definition}[Degeneracy Relations for the Morse-Bott-Smale Chain Complex]\label{degenerate}
Let $\sigma_P, \sigma_Q \in S_p^\infty(B_i)$ be singular $C_p$-spaces in $B_i$ and let
$\partial Q = \sum_j n_j Q_j \in S_{p-1}$.  For any map $\alpha:P \rightarrow Q$,
let $\partial_0 \sigma_Q \circ \alpha$ denote the formal sum $(-1)^{p+i} \sum_j n_j
(\sigma_Q\circ\alpha) |_{\alpha^{-1}(Q_j)}$. Define the subgroup 
$D_p^\infty(B_i)\subseteq S_p^\infty(B_i)$ of \textbf{degenerate singular 
topological chains} to be the subgroup generated by the following elements. 
\begin{enumerate}
\item If $\alpha$ is an orientation preserving homeomorphism such that
      $\sigma_Q \circ \alpha = \sigma_P$ and $\partial_0 \sigma_Q \circ 
      \alpha = \partial_0 \sigma_P$, then $\sigma_P - \sigma_Q \in D_p^\infty(B_i)$.
\item If $P$ is a face of $I^N$ and $\sigma_P$ does not depend on some free 
      coordinate of $P$, then $\sigma_P \in D_p^\infty(B_i)$ and
      $\partial_j(\sigma_P) \in D_{p+j-1}^\infty(B_{i-j})$ for all $j=1,\ldots ,m$.  
\item If $P$ and $Q$ are connected components of some fibered products and
      $\alpha$ is an orientation reversing map such that
      $\sigma_Q \circ \alpha = \sigma_P$ and $\partial_0 \sigma_Q \circ 
      \alpha = \partial_0 \sigma_P$, then $\sigma_P + \sigma_Q \in D_p^\infty(B_i)$.
\item If $Q$ is a face of $I^N$ and $R$ is a connected component of a fibered product
      $$
      Q \times_{B_{i_1}} \overline{\mathcal{M}}(B_{i_1},B_{i_2}) \times_{B_{i_2}}
      \overline{\mathcal{M}}(B_{i_2},B_{i_3}) \times_{B_{i_3}} \cdots 
      \times_{B_{i_{n-1}}} \overline{\mathcal{M}}(B_{i_{n-1}},B_{i_n})
      $$
      such that $\text{deg }R > \text{dim }B_{i_n}$, then $\sigma_R \in D_r^\infty(B_{i_n})$
      and $\partial_j(\sigma_R) \in D_{r+j-1}^\infty(B_{{i_n}-j})$ for all
      $j = 0,\ldots ,m$.
\item If $\sum_\alpha n_\alpha \sigma_\alpha \in S_\ast(R)$ 
      is a smooth singular chain in a connected component $R$ of a 
      fibered product (as in (4)) that represents the fundamental class of $R$ and 
      $$
      (-1)^{r + i_n} \partial_0 \sigma_R - \sum_\alpha n_\alpha \partial (\sigma_R \circ \sigma_\alpha)
      $$
      is in the group generated by the elements satisfying one of the above conditions,
      then
      $$
      \sigma_R - \sum_\alpha n_\alpha (\sigma_R \circ \sigma_\alpha) \in D_r^\infty(B_{i_n})
      $$
      and
      $$
      \partial_j \left(\sigma_R - \sum_\alpha n_\alpha (\sigma_R \circ \sigma_\alpha) \right)
      \in D_{r+j-1}^\infty(B_{{i_n}-j})
      $$
      for all $j = 1,\ldots ,m$.
\end{enumerate}
\end{definition}

\smallskip\noindent
Note: Condition 3 does not apply to the $p$-faces of $I^N$.
If we allowed condition 3 to apply to the $p$-faces of $I^N$, then
we could combine conditions 1 and 3.  However, the resulting chain 
complex would then reduce to a quotient of the chain complex of smooth singular
cubes when $f:M \rightarrow \mathbb{R}$ is a constant function
(see Example \ref{constantfunction}).

\begin{lemma}\label{boundarydegenerate}
For any $i,j=0,\ldots , m$ there is an induced homomorphism
$$
\partial_j: S_{p}^\infty(B_i)/D_p^\infty(B_i) \rightarrow
S_{p+j-1}^\infty(B_{i-j})/D_{p+j-1}^\infty(B_{i-j}).
$$
\end{lemma}

\smallskip\noindent
Proof:
Let $\sigma_P:P \rightarrow B_i$ be a singular $C_p$-space in $S_p^\infty(B_i)$
such that $\sigma_Q \circ \alpha = \sigma_P$ and $\partial_0 \sigma_Q \circ \alpha =
\partial_0 \sigma_P$ for some singular $C_p$-space $\sigma_Q:Q \rightarrow B_i$
in $S_p^\infty(B_i)$ and some orientation preserving homeomorphism $\alpha:P \rightarrow Q$.
By Proposition \ref{complex}, $\partial^2_0\sigma_Q = \partial^2_0 \sigma_P = 0$ and
hence, $\partial_0(\partial_0 \sigma_Q) \circ \alpha = \partial_0(\partial_0 \sigma_P$).
Thus, $\partial_0(\sigma_P - \sigma_Q) = \partial_0 \sigma_P - \partial_0\sigma_Q 
\in D_{p-1}^\infty(B_i)$.
Now assume that  $1 \leq j \leq i$ and recall that 
$\partial_j(\sigma_P)$ is the map
$$
\partial_j(\sigma_P): P \times_{B_i}\overline{\mathcal{M}}(B_i,B_{i-j})
\stackrel{\pi_2}{\longrightarrow}\overline{\mathcal{M}}(B_i,B_{i-j})
\stackrel{\partial_+}{\longrightarrow} B_{i-j}
$$
where the fibered product is taken over $\sigma_P:P \rightarrow B_i$
and $\partial_-:\overline{\mathcal{M}}(B_i,B_{i-j}) \rightarrow B_i$.
The map $\beta:P \times_{B_i}\overline{\mathcal{M}}(B_i,B_{i-j}) \rightarrow
Q \times_{B_i}\overline{\mathcal{M}}(B_i,B_{i-j})$ defined by
$\beta(p,\gamma) = (\alpha(p),\gamma)$ is an orientation preserving
homeomorphism that satisfies $\partial_j(\sigma_Q) \circ \beta = 
\partial_j(\sigma_P)$ since $\sigma_Q \circ \alpha = \sigma_P$. By 
Definition \ref{fiberedboundary},
$$
\partial(P \times_{B_i} \overline{\mathcal{M}}(B_i,B_{i-j})) =
\partial P \times_{B_i} \overline{\mathcal{M}}(B_i,B_{i-j}) +
(-1)^{p+b_i} P \times_{B_i} \partial\overline{\mathcal{M}}(B_i,B_{i-j})
$$
and
$$
\partial(Q \times_{B_i} \overline{\mathcal{M}}(B_i,B_{i-j})) =
\partial Q \times_{B_i} \overline{\mathcal{M}}(B_i,B_{i-j}) +
(-1)^{p+b_i} Q \times_{B_i} \partial\overline{\mathcal{M}}(B_i,B_{i-j})
$$
The assumption $\partial_0 \sigma_Q \circ \alpha = \partial_0
\sigma_P$ implies that $\partial_0 (\partial_j(\sigma_Q)) \circ \beta
=\partial_0(\partial_j(\sigma_P))$, and thus $\partial_j(\sigma_P - \sigma_Q)
= \partial_j \sigma_P - \partial_j\sigma_Q \in D_{p+j-1}^\infty(B_{i-j})$
by the first condition for degeneracy.  This completes the proof for
the first case.  The proof for the third case is similar.

Now assume that $\sigma_P:P \rightarrow B_i$ is a singular $C_p$-space
in $S_p^\infty(B_i)$ where $P$ is a face of $I^N$ and $\sigma_P$ does 
not depend on some free coordinate $x_j$. Then in the sum for $\partial_0 
\sigma_P$ the term $\sigma_P|_{x_j = 0} - \sigma_P|_{x_j = 1} \in 
D_{p-1}^\infty(B_i)$ because it satisfies the first condition for 
degeneracy, and the rest of the terms are all independent of the coordinate 
corresponding to $x_j$ and thus satisfy the second condition for degeneracy. 
Thus, $\partial_0 \sigma_P \in D_{p-1}^\infty(B_i)$. If $j>0$, then the 
second condition for degeneracy implies that $\partial_j(\sigma_P) \in 
D_{p+j-1}^\infty(B_{i-j})$. This completes the proof for the second case 
of Definition \ref{degenerate}.  The fourth and fifth cases follow 
immediately from the definition.

\begin{flushright}
$\Box$
\end{flushright}

\begin{definition}\label{MBScomplex}
Define
$$
C_p(B_i) = S_p^\infty(B_i) / D_p^\infty(B_i)
$$
to be the group of \textbf{non-degenerate} smooth singular topological chains in
$S_p^\infty(B_i)$. The group $C_k(f)$ of \textbf{$k$-chains} in the Morse-Bott
chain complex of $f$ is defined to be the group of non-degenerate smooth singular
topological chains of Morse-Bott degree $k$, i.e.
$$
C_k(f) = \bigoplus_{i=0}^m C_{k-i}(B_i) = \bigoplus_{i=0}^m
S_{k-i}^\infty(B_i)/D_{k-i}^\infty(B_i).
$$
The boundary operator in the Morse-Bott-Smale chain complex
$$
\partial:\bigoplus_{i=0}^m S_{k-i}^\infty (B_i)/D_{k-i}^\infty(B_i) \longrightarrow
\bigoplus_{i=0}^m S_{k-1-i}^\infty(B_i)/D_{k-1-i}^\infty(B_i)
$$
is defined to be $\partial = \oplus_{j=0}^m \partial_j$
where $\partial_j$ is the induced homomorphism from Lemma 
\ref{boundarydegenerate}.
\end{definition}

\noindent
When $m=2$ this homomorphism can be pictured as follows.

$$
\xymatrix{
S_0^\infty(B_2)/D_0^\infty(B_2) \ar[r]^{\partial_0} \ar[dr]^{\partial_1} \ar[ddr]|(.33){\partial_2} |!{[d];[dr]}\hole & 0 & & \\
S_1^\infty(B_1)/D_1^\infty(B_1) \ar@{}[u]|{\oplus} \ar[r]^(.55){\partial_0} \ar[dr]^{\partial_1} & S_0^\infty(B_1)/D_0^\infty(B_1)
   \ar@{}[u]|{\oplus} \ar[r]^{\partial_0} \ar[dr]^{\partial_1} & 0 & \\
S_2^\infty(B_0)/D_2^\infty(B_0) \ar@{}[u]|{\oplus} \ar[r]^{\partial_0} & S_1^\infty(B_0)/D_1^\infty(B_0) \ar@{}[u]|{\oplus} \ar[r]^{\partial_0} & 
   S_0^\infty(B_0)/D_0^\infty(B_0) \ar@{}[u]|{\oplus}  \ar[r]^-{\partial_0} & 0\\
C_2(f) \ar@{}[u]|{\|} \ar[r]^{\partial} & C_1(f) \ar@{}[u]|{\|} \ar[r]^{\partial} & C_0(f) \ar@{}[u]|{\|} \ar[r]^-{\partial} & 0
}
$$

\begin{corollary}\label{MBisChain}
The pair $(C_\ast(f),\partial)$ is a chain complex, 
i.e. $\partial \circ \partial = 0$. 
\end{corollary}

\smallskip\noindent
Proof:
This follows from Proposition \ref{complex} and Lemma
\ref{boundarydegenerate}.

\begin{flushright}
$\Box$
\end{flushright}


\subsection{Computing Morse-Bott homology}\label{examples}

\noindent
\begin{example}[A constant function]\label{constantfunction}
\end{example}

Let $f:M \rightarrow \mathbb{R}$ be a constant function on a compact
oriented smooth Riemannian manifold $(M,g)$, and let $N > \text{dim }M$. 
Since $f$ is constant there is only one critical submanifold
$B_0 = M$, and for any $p \geq 0$ the set $C_p$ consists of the faces
of $I^N$ of dimension $p$. For any $k \geq 0$ the group $\tilde{C}_k(f) = S_k^\infty(B_0)$
is the subgroup of smooth chains in the group $S_k(B)$ of continuous singular
$N$-cube chains of degree $k$, $\partial_0 = (-1)^k \partial$, and $\partial_j$
is trivial for all $j >0$.  

Conditions 1 and 2 for degeneracy in Definition \ref{degenerate} 
agree with the two conditions for degeneracy in Definition 
\ref{cubedegenerate}, and the other conditions are vacuous since 
$B_i = \emptyset$ for all $i > 0$. Thus, $D_k^\infty(B_0)$ is the 
subgroup smooth chains in the group $D_k(B_0)$ of continuous degenerate 
singular $N$-cube chains, and the Morse-Bott-Smale chain complex $(C_\ast(f),\partial)$,
with its boundary operator multiplied by $(-1)^k$, is the subchain complex 
of the chain complex $(S_\ast(B)/D_\ast(B),\partial_\ast)$ from 
Theorem \ref{cubehomology} consisting of nondegenerate smooth 
singular $N$-cube chains. The proof of Theorem \ref{cubehomology} 
carries over verbatim to the smooth case, and standard arguments 
can be used to show that the homology of the \textbf{smooth} singular
cube chain complex of a smooth manifold $M$ is isomorphic to the 
homology of the \textbf{continuous} singular cube chain complex, cf. 
Chapter 5 of \cite{WarFou}. Thus, the homology of the Morse-Bott-Smale
chain complex $(C_\ast(f),\partial)$ is isomorphic to the singular 
homology $H_\ast(M;\mathbb{Z})$. When $m=2$ this chain complex can 
be pictured as follows.

$$
\xymatrix{
0 \ar[r]^{\partial_0} \ar[dr]^{\partial_1} \ar[ddr]|(.33){\partial_2} |!{[d];[dr]}\hole & 0 & & \\
0 \ar@{}[u]|{\oplus} \ar[r]^(.55){\partial_0} \ar[dr]^{\partial_1} & 0
   \ar@{}[u]|{\oplus} \ar[r]^{\partial_0} \ar[dr]^{\partial_1} & 0 & \\
S_2^\infty(B_0)/D_2^\infty(B_0) \ar@{}[u]|{\oplus} \ar[r]^{\partial_0} & S_1^\infty(B_0)/D_1^\infty(B_0) \ar@{}[u]|{\oplus} \ar[r]^{\partial_0} & 
   S_0^\infty(B_0)/D_0^\infty(B_0) \ar@{}[u]|{\oplus}  \ar[r]^-{\partial_0} & 0\\
C_2(f) \ar@{}[u]|{\|} \ar[r]^{\partial} & C_1(f) \ar@{}[u]|{\|} \ar[r]^{\partial} & C_0(f) \ar@{}[u]|{\|} \ar[r]^-{\partial} & 0
}
$$

\noindent
\begin{example}[A Morse-Smale function]\label{MSfunction}
\end{example}

Let $f:M \rightarrow \mathbb{R}$ be a Morse-Smale function on a compact
oriented smooth Riemannian manifold $(M,g)$, and let $N > \text{dim }M$. 
Since the critical points of $f$ are isolated, $\text{dim }B_i = 0$
for all $i$ and conditions 2 and 4 for degeneracy in Definition \ref{degenerate}
imply that $D_p^\infty(B_i) = S_p^\infty(B_i)$ for all $i \geq 0$ 
and $p > 0$. This implies that $\partial_j$ is trivial unless $j=1$
(see Lemma \ref{partialj}). 
 
Moreover, the first condition for degeneracy in Definition \ref{degenerate}
implies that $C_k(f) = S_0(B_k)/D_0(B_k)$ is isomorphic to the free 
abelian group generated by the critical points of index $k$ for all 
$k=0,\ldots ,m$.  If $q \in B_k$ is a critical point of index $k$, 
then the map $\partial_1$ applied to $q$:
$$
\partial_1(q):q \times_{B_k} \overline{\mathcal{M}}(B_k,B_{k-1}) 
\stackrel{\pi_2}{\longrightarrow} \overline{\mathcal{M}}(B_k,B_{k-1})
\stackrel{\partial_+}{\longrightarrow} B_{k-1}
$$
counts the number of gradient flow lines (with sign) between
$q$ and the critical points in $B_{k-1}$, i.e.
$$
\partial_1(q) = \sum_{p \in B_{k-1}} n(q,p)p \in C_{k-1}(f)
$$
where $n(q,p)$ is the number of gradient flow lines between 
$q$ and $p$ counted with the signs coming from the orientation
of the zero dimensional manifold $\overline{\mathcal{M}}(B_k,B_{k-1})$.
Hence, $(C_\ast(f),\partial)$ is the Morse-Smale-Witten
chain complex, and the homology of $(C_\ast(f),\partial)$ is 
isomorphic to the singular homology $H_\ast(M;\mathbb{Z})$ by
the Morse Homology Theorem. When $m=2$ this chain complex can 
be pictured as follows.

$$
\xymatrix{
S_0^\infty(B_2)/D_0^\infty(B_2) \ar[r]^{\partial_0} \ar[dr]^{\partial_1} \ar[ddr]|(.33){\partial_2} |!{[d];[dr]}\hole & 0 & & \\
0 \ar@{}[u]|{\oplus} \ar[r]^(.55){\partial_0} \ar[dr]^{\partial_1} & S_0^\infty(B_1)/D_0^\infty(B_1)
   \ar@{}[u]|{\oplus} \ar[r]^{\partial_0} \ar[dr]^{\partial_1} & 0 & \\
0 \ar@{}[u]|{\oplus} \ar[r]^{\partial_0} & 0 \ar@{}[u]|{\oplus} \ar[r]^{\partial_0} & 
   S_0^\infty(B_0)/D_0^\infty(B_0) \ar@{}[u]|{\oplus}  \ar[r]^-{\partial_0} & 0\\
C_2(f) \ar@{}[u]|{\|} \ar[r]^{\partial} & C_1(f) \ar@{}[u]|{\|} \ar[r]^{\partial} & C_0(f) \ar@{}[u]|{\|} \ar[r]^-{\partial} & 0
}
$$

\noindent
\begin{example}[A Morse-Bott-Smale function on $S^2$]\label{MSBsphere1}
\end{example}

Consider $M = S^2 = \{(x,y,z) \in \mathbb{R}^3 |\ x^2 + y^2 + z^2 = 1 \}$,
and let $f(x,y,z) = z^2$. Then $B_0 \approx S^1$, $B_1 = \emptyset$,
and $B_2 = \{n,s\}$ where $n=(0,0,1)$ and $s=(0,0,-1)$.

\begin{center}
\includegraphics{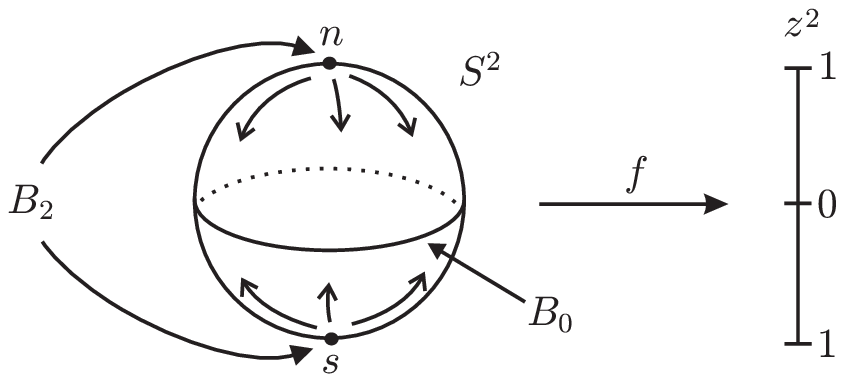}
\end{center}

\noindent
Note that
$$
S_0^\infty(B_2)/D_0^\infty(B_2) \approx <n,s> \approx \mathbb{Z} \oplus \mathbb{Z}
$$
by the first two conditions for degeneracy, and $S_p^\infty(B_2)/D_p^\infty(B_2) = 0$
for $p > 0$ by the second condition for degeneracy in Definition \ref{degenerate}.
Thus, the Morse-Bott-Smale chain complex can be pictured as follows.

$$
\xymatrix{
<n,s> \ar[r]^{\partial_0} \ar[dr]^{\partial_1} \ar[ddr]|(.33){\partial_2} |!{[d];[dr]}\hole & 0 & & \\
0 \ar@{}[u]|{\oplus} \ar[r]^(.55){\partial_0} \ar[dr]^{\partial_1} & 0 \ar@{}[u]|{\oplus} \ar[r]^{\partial_0} \ar[dr]^{\partial_1} & 0 & \\
S_2^\infty(B_0)/D_2^\infty(B_0) \ar@{}[u]|{\oplus} \ar[r]^{\partial_0} & S_1^\infty(B_0)/D_1^\infty(B_0) \ar@{}[u]|{\oplus} \ar[r]^{\partial_0} & 
   S_0^\infty(B_0)/D_0^\infty(B_0) \ar@{}[u]|{\oplus}  \ar[r]^-{\partial_0} & 0\\
C_2(f) \ar@{}[u]|{\|} \ar[r]^{\partial} & C_1(f) \ar@{}[u]|{\|} \ar[r]^{\partial} & C_0(f) \ar@{}[u]|{\|} \ar[r]^-{\partial} & 0
}
$$

\noindent
The group $S_k^\infty(B_0)/D_k^\infty(B_0)$ is non-trivial for all $k\leq N$, but
$$
H_k(C_\ast(f),\partial) = 0 \quad \text{if} \quad k > 2
$$
and $\partial_0:S_3^\infty(B_0)/D_3^\infty(B_0) \rightarrow S_2^\infty(B_0)/D_2^\infty(B_0)$
maps onto the kernel of the boundary operator 
$\partial_0:S_2^\infty(B_0)/D_2^\infty(B_0) \rightarrow S_1^\infty(B_0)/D_1^\infty(B_0)$
because the homology of the bottom row in the above diagram is the smooth integral
singular homology of $B_0 \approx S^1$.

The moduli space $\overline{\mathcal{M}}(B_2,B_0)$ is a disjoint union of two copies
of $S^1$ with opposite orientations. This moduli space can be viewed as a subset of
the manifold $S^2$ since $\overline{\mathcal{M}}(B_2,B_0) = \mathcal{M}(B_2,B_0)$.

\begin{center}
\includegraphics{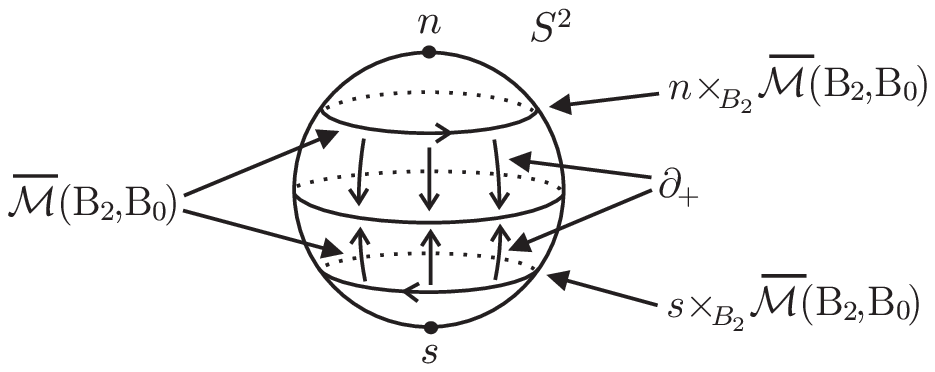}
\end{center}

\noindent
Thus, there is an orientation reversing map $\alpha:n \times_n
\overline{\mathcal{M}}(B_2,B_0) \rightarrow s \times_s \overline{\mathcal{M}}(B_2,B_0)$
such that $\partial_2(n) \circ \alpha = \partial_2(s)$.  Since
$\partial_0(\partial_2(n)) = \partial_0(\partial_2(s)) = 0$, the third condition
for degeneracy from Definition \ref{degenerate} implies that
$$
\partial_2(n + s) = \partial_2(n) + \partial_2(s) = 0 \in S_1(B_0)/D_1(B_0),
$$
and the fifth condition for degeneracy implies that $\partial_2$ maps either 
$n$ or $s$ onto a representative of the generator of 
$$
\frac{\text{ker }\partial_0:S_1^\infty(B_0)/D_1^\infty(B_0) \rightarrow S_0^\infty(B_0)/D_0^\infty(B_0)}
{\text{im }\partial_0:S_2^\infty(B_0)/D_2^\infty(B_0) \rightarrow S_1^\infty(B_0)/D_1^\infty(B_0)}
\approx H_1(S^1;\mathbb{Z}) \approx \mathbb{Z}
$$
depending on the orientation chosen for $B_0$. Therefore,
$$
H_k(C_\ast(f),\partial) = \left\{
\begin{array}{ll}
\mathbb{Z} & \text{if } k = 0,2\\
0          & \text{otherwise.}
\end{array}\right.
$$

\noindent
\begin{example}[Another Morse-Bott-Smale function on $S^2$]\label{MSBsphere2}
\end{example}

Consider $M = S^2 = \{(x,y,z) \in \mathbb{R}^3 |\ x^2 + y^2 + z^2 = 1 \}$,
and let $f(x,y,z) = -z^2$. Then $B_0 = \{n,s\}$ where $n=(0,0,1)$ and 
$s=(0,0,-1)$, $B_1 \approx S^1$, and $B_2 = \emptyset$.

\begin{center}
\includegraphics{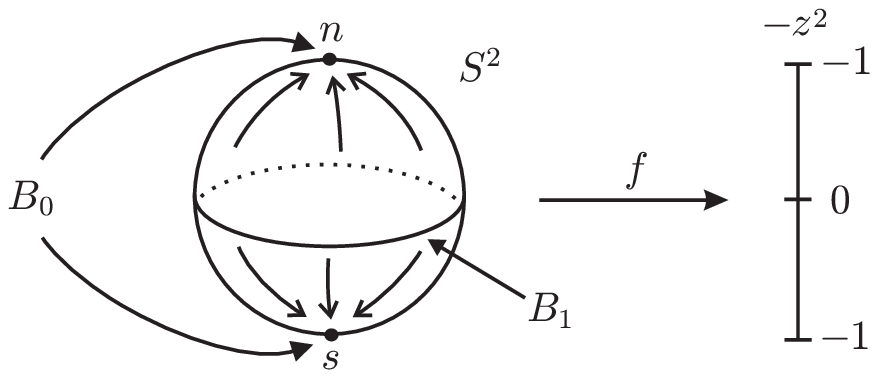}
\end{center}

\noindent
We have 
$$
S_0^\infty(B_0)/D_0^\infty(B_0) \approx <n,s> \approx \mathbb{Z} \oplus \mathbb{Z}
$$
by the first two conditions for degeneracy, and $S_p^\infty(B_0)/D_p^\infty(B_0) = 0$
for $p > 0$ by the second and fourth conditions for degeneracy in 
Definition \ref{degenerate}. Thus, the Morse-Bott-Smale chain complex can be
pictured as follows.

$$
\xymatrix{
0 \ar[r]^{\partial_0} \ar[dr]^{\partial_1} \ar[ddr]|(.33){\partial_2} |!{[d];[dr]}\hole & 0 & & \\
S_1^\infty(B_1)/D_1^\infty(B_1) \ar@{}[u]|{\oplus} \ar[r]^(.55){\partial_0} \ar[dr]^{\partial_1} & S_0^\infty(B_1)/D_0^\infty(B_1)
   \ar@{}[u]|{\oplus} \ar[r]^{\partial_0} \ar[dr]^{\partial_1} & 0 & \\
0 \ar@{}[u]|{\oplus} \ar[r]^{\partial_0} & 0 \ar@{}[u]|{\oplus} \ar[r]^{\partial_0} & 
   <n,s> \ar@{}[u]|{\oplus}  \ar[r]^-{\partial_0} & 0\\
C_2(f) \ar@{}[u]|{\|} \ar[r]^{\partial} & C_1(f) \ar@{}[u]|{\|} \ar[r]^{\partial} & C_0(f) \ar@{}[u]|{\|} \ar[r]^-{\partial} & 0
}
$$

\noindent
The middle row computes the smooth integral singular homology of $B_1 \approx S^1$
when $k>1$, and hence $H_2(C_\ast(f),\partial) = \mathbb{Z}$.  
However, for $k=1$ the kernel of $\partial = \partial_0 \oplus \partial_1:
C_1(f) \rightarrow C_0(f)$ is trivial because $\partial_1(\sigma_P)$ is
non-zero for any map $\sigma_P:P \rightarrow B_1$ from a point $P\in C_0$
to $B_1$. In fact, the fibered product $P \times_{B_1} \overline{\mathcal{M}}(B_1,B_0)$
consists of two points, and
$$
\partial_1(\sigma_P):P \times_{B_1} \overline{\mathcal{M}}(B_1,B_0)
\stackrel{\pi_2}{\longrightarrow} \overline{\mathcal{M}}(B_1,B_0)
\stackrel{\partial_+}{\longrightarrow} \{n,s\}
$$ 
represents either $n-s \in C_0(f)$ or $s-n \in C_0(f)$, depending
on the orientation chosen for $B_1$. (If we fix an orientation on
$B_1$, then the sign is determined by Definition \ref{fiberedorientation}).
Therefore,
$$
H_k(C_\ast(f),\partial) = \left\{
\begin{array}{ll}
\mathbb{Z} & \text{if } k = 0,2\\
0          & \text{otherwise.}
\end{array}\right.
$$


\subsection{Manifolds with Corners}\label{manicorners}

In this subsection we recall some facts about manifolds with corners,
and we prove Lemmas \ref{basiscorners} and \ref{fiberedmanifold}. 
Note that the definition of the Morse-Bott-Smale chain complex does 
not rely on the smooth manifold with corners structure on the objects 
in $C_p$. All that is required (in addition to the abstract 
topological chain structure) is the notion of an orientation on these 
objects and some finiteness conditions. In fact, to define the 
Morse-Bott-Smale chain complex it would be sufficient to show that the 
objects in $C_p$ are oriented topological manifolds with boundary 
and for any $P \in C_p$ 
$$
P \times_{B_i} \overline{\mathcal{M}}(B_i,B_{i-j})
$$
has finitely many components.

However, the following example (adapted from \cite{BauTop}) shows that some 
additional structure is needed to guarantee that these fibered
products have finitely many components.
\begin{example}\label{badfibered}
Let $f: [0,1] \rightarrow [0,1] \times [-1,1]$ be given by 
$$
f(t) = \left\{\begin{array}{ll}
(t, e^{-1/t^2} \sin(\pi/t) ) & \text{if } t \neq 0 \\
(0,0) & \text{if }t = 0 
\end{array}\right.
$$
and $g:[0,1] \times [0,1] \rightarrow [0,1] \times [-1,1]$ be given
by $g(x,y) = (x,0)$. Then $f$ and $g$ are smooth maps between finite
dimensional compact oriented topological manifolds with boundary 
whose fibered product
$$
[0,1] \times _{(f,g)} [0,1] \times [0,1] = \{(t,t,0) \in [0,1] \times [0,1] \times [0,1] 
|\ t = 0,1,1/2,1/3, \ldots  \}
$$
has infinitely many components. 
\end{example}

\noindent
Note that this example also shows that the fibered product of
two finite CW-complexes might not be a CW-complex, and the
fibered product of two finite simplicial complexes might
not be a finite simplicial complex. So, a category with some
additional structure is needed in order to prove that the fibered products
of interest for Morse-Bott homology have finitely many components.
An analysis of the above example shows that the failure of this fibered 
product to have finitely many components is due to a lack of transversality.
Thus, the category of compact smooth manifolds with corners is a good 
candidate for proving the required finiteness conditions. Other choices
are possible, including the category of compact manifolds with corners
of class $C^s$ where $s \geq 1$.

\smallskip
Recall that a manifold with corners of dimension $m$ is a Hausdorff
topological space $X$ such that every point $x \in X$ has a neighborhood
$U_x\subseteq X$ that is locally  homeomorphic via a chart
$\psi_x:U_x \rightarrow \mathbb{R}^m_{(r)}$ to some open subset of 
$\mathbb{R}^m_{(r)} = \mathbb{R}^r \times [0,\infty)^{m-r}$.
A manifold with corners is said to be of class $C^s$ if and only if
for any two charts $\psi_x:U_x \rightarrow \mathbb{R}^m_{(r)}$ and 
$\psi_y:U_y \rightarrow \mathbb{R}^m_{(r)}$ the composition $\psi_x 
\circ \psi_y^{-1}: \psi_y(U_x \cap U_y) \rightarrow \psi_x(U_x \cap U_y)$
extends to a map of class $C^s$ between open sets of $\mathbb{R}^m$.
We will assume that $s \geq 1$, and if $s=\infty$, then we will call
$X$ a smooth manifold with corners.

If $\psi_x:U_x \rightarrow \mathbb{R}^m_{(r)}$ is a local chart on $X$
such that $\psi_x(x) = 0 \in \mathbb{R}^m_{(r)}$, then the number
$r$ does not depend on the choice of the chart $\psi_x$ and is
called the \textbf{index} of $x$ in $X$, denoted $\mbox{ind}(X,x)$.
A manifold with corners $X$ is stratified by manifolds (without boundary)
$X_0, X_1, \ldots ,X_m$ defined by
$$
X_r = \{x \in X|\ \text{ind}(X,x) = r\}
$$
for all $r=0,\ldots,m$. A connected component of $X_r$ is called an
$r$-\textbf{stratum} of $X$, and a subset of $X$ that is an $r$-stratum
for some $r$ is called a \textbf{stratum} of $X$.  The set $X_m$ is 
the \textbf{interior} of $X$, and the set $X - X_m$ is called the 
\textbf{boundary} if $X$. For more details see  \cite{LatGra}, 
\cite{MarDif}, or \cite{NieTra}.

The strata of the compact smooth manifold with corners 
$\overline{\mathcal{M}}(B_i,B_{i-j})$ can be described as follows. 
To simplify the notation assume that the components of $B_i$ are 
all of dimension $b_i$, which implies that 
$\overline{\mathcal{M}}(B_i,B_{i-j})$ has dimension $m = b_i+j-1$. 
The strata are then the connected components of the sets
$$
X_{m-r} = \bigcup_{0 <j_1 < \cdots <j_r < j} \mathcal{M}(B_i,B_{i-j_1})
\times_{B_{i-j_1}} \mathcal{M}(B_{i-j_1},B_{i-j_2}) \times_{B_{i-j_2}} \cdots
\times_{B_{i-j_r}} \mathcal{M}(B_{i-j_r},B_{i-j})
$$
where $r=0,\ldots ,j-1$, and the union is taken over all increasing sequences
of integers $j_1<j_2 <\cdots < j_r$ between $1$ and $j-1$. Note that
$X_m = \mathcal{M}(B_i,B_{i-j})$, and if $1 \leq r \leq j-1$, then the set $X_{m-r}$ 
consists of those piecewise gradient flow lines that pass through exactly 
$r$ intermediate critical submanifolds.

\smallskip
We will use the following theorem (see Theorem 3 of \cite{NieTra}).

\begin{theorem}\label{preimagecorners}
Let $X$ and $Y$ be $C^s$ manifolds with corners, where $s \geq 1$.
Let $A \subseteq Y$ be a $C^s$ submanifold with corners, and 
$f:X \rightarrow Y$ a local $C^s$ map, which preserves local facets 
relatively to $A$ and intersects $A$ transversally and stratum transversally.
Then either $f^{-1}(A) = \emptyset$, or
\begin{enumerate}
\item $f^{-1}(A)$ is a $C^s$ submanifold with corners of $X$, and
\item $\text{dim }X - \text{dim }f^{-1}(A) = \text{dim }Y - \text{dim }A$, and
\item $\text{ind}(X,x) - \text{ind}(f^{-1}(A),x) = \text{ind}(Y,f(x)) - \text{ind}(A,f(x))$
      for all $x \in f^{-1}(A)$.
\end{enumerate}
\end{theorem}

\noindent
For a definition of ``local facets" see \cite{NieTra}. We will not 
define local facets here because we will only apply this theorem 
in the case where $Y$ is a manifold without boundary, and the local
facets condition is always satisfied in that case. The assumption that 
$f$ intersects $A$ \textbf{stratum transversally} means that for any 
$x \in f^{-1}(A)$ we have
$$
df_x(\hat{T}_xX) + \hat{T}_yA = \hat{T}_yY
$$
where $y = f(x)$ and $\hat{T}_xX$ denotes the tangent space of the 
stratum containing $x \in X$. Similarly, we will say that 
a map $f:X \rightarrow Y$ is a \textbf{stratum submersion} at $x \in X$
if and only if $df_x$ maps $\hat{T}_xX$ onto $\hat{T}_yY$ where $y=f(x)$.
Note that if $f$ is a stratum submersion at $x \in X$ and $A \subseteq Y$
is any submanifold with corners containing $y$, then $f$ intersects 
$A$ stratum transversally.

\begin{lemma}\label{beginsubmersioninterior}
For any two connected critical submanifolds $B$ and $B'$, the beginning
point map
$$
\partial_-:\mathcal{M}(B,B') \rightarrow B
$$
is a submersion.
\end{lemma}

\smallskip\noindent
Proof: Let $B$ be a connected critical submanifold of a
Morse-Bott-Smale function $f:M \rightarrow \mathbb{R}$
of dimension $b$ and index $\lambda_B$. By Theorem 
\ref{stablemanifold} the beginning point map
$$
\partial_-:W^u(B) \rightarrow B
$$
is a locally trivial fiber bundle, and if $B'$ is any other critical 
submanifold, then either $W(B,B') = \emptyset$ (in which case the 
lemma holds trivially) or $W(B,B')$ is a manifold of dimension
$\lambda_B - \lambda_{B'} + b$ (by Lemma \ref{dim}).

Let $x \in B$ and let $U\subset B$ be a coordinate neighborhood
of $x$ such that the bundle $\partial_-:W^u(B) \rightarrow B$
is trivial over $U$. We can extend the coordinate chart $U$
on $B$ to a coordinate chart $V$ on $M$ so that in $V$ we 
have local coordinates
$$
\mathbb{R}^b \times \mathbb{R}^{\lambda_B} \times \mathbb{R}^{m-b-\lambda_B}
$$
where $\mathbb{R}^b$ represents the coordinates on $B$ and $\mathbb{R}^{\lambda_B}$
represents the coordinates along the fiber of $W^u(B)$. In these
local coordinates the beginning point map $\partial_-:W^u(U) \rightarrow U$
is given by the projection map
$$
\pi: \mathbb{R}^b \times \mathbb{R}^{\lambda_B} \times \vec{0}
\rightarrow \mathbb{R}^b \times \vec{0} \times \vec{0}.
$$
 
Let $\tilde{x} \in V$ be a point in the fiber above $x$, i.e.
$\partial_-(\tilde{x}) = x$, such that $\tilde{x} \in W^s(B')$.
In the local coordinates we have
$T_{\tilde{x}} W^u(x) = \vec{0} \times \mathbb{R}^{\lambda_B} \times \vec{0}$,
and the Morse-Bott-Smale transversality condition 
$$
(\vec{0} \times \mathbb{R}^{\lambda_B} \times \vec{0}) + T_{\tilde{x}} W^s(B') = \mathbb{R}^m
\approx T_{\tilde{x}} M
$$
implies that we can find a smooth section $s:U' \rightarrow V$
of the bundle $\partial_-:W^u(B) \rightarrow B$ on an open set 
$U' \subseteq U$ containing $x$ such that $s(x) = \tilde{x}$ 
and $s(z) \in W^s(B')$ for all $z \in U'$.
The image of the section $s:U' \rightarrow V$ is in
$W(B,B') = W^u(B) \cap W^s(B')$ by construction, and 
diffeomorphic to the open set $U' \subseteq B$ via
the beginning point map $\partial_-:W^u(U') \rightarrow U'$.
This shows that the beginning point map
$$
\partial_-:W(B,B') \rightarrow B
$$
is a submersion at $\tilde{x}$, and using the diffeomorphisms
determined by the gradient flow we see that the beginning
point map is a submersion at every point $\tilde{x} \in W(B,B')$.

If we pick an appropriate level set $f^{-1}(a)$ between $B$ and 
$B'$, then we can identify
$$
\mathcal{M}(B,B') = W(B,B') \cap f^{-1}(a).
$$
This shows that 
$$
\partial_-:\mathcal{M}(B,B') \rightarrow B
$$
is a submersion because the level sets are transverse
to the gradient flow.

\begin{flushright}
$\Box$
\end{flushright}

\begin{corollary}\label{beginsubmersion}
For any two connected critical submanifolds $B$ and $B'$, the beginning
point map
$$
\partial_-:\overline{\mathcal{M}}(B,B') \rightarrow B
$$
is a submersion and a stratum submersion.
\end{corollary}

\smallskip\noindent
Proof:
Since $B$ is a manifold without boundary, $\hat{T}_yB = T_yB$ 
for all $y \in B$.  Hence, it suffices to show that 
$\partial_-:\overline{\mathcal{M}}(B,B') \rightarrow B$
is a stratum submersion. The boundary of 
$\overline{\mathcal{M}}(B,B')$ is stratified by fibered products 
of moduli spaces of gradient flow lines of the form
$$
\mathcal{M}(B,B_1) \times_{B_1} \mathcal{M}(B_1,B_2) \times_{B_2}
\cdots \times_{B_n} \mathcal{M}(B_n,B')
$$
for some intermediate critical submanifolds $B_1,\ldots ,B_n$.
To see that these strata are smooth manifolds without
boundary assume for the purpose of induction that
$$
\mathcal{M}(B,B_1) \times_{B_1} \mathcal{M}(B_1,B_2) \times_{B_2}
\cdots \times_{B_{n-1}} \mathcal{M}(B_{n-1},B_n)
$$
is a smooth manifold without boundary. By Lemma \ref{beginsubmersioninterior}
the map $\partial_-:\mathcal{M}(B_n,B') \rightarrow B_n$
is a submersion, and hence
$$
\mathcal{M}(B,B_1) \times_{B_1} \mathcal{M}(B_1,B_2) \times_{B_2}
\cdots \times_{B_n} \mathcal{M}(B_n,B')
$$
is a smooth manifold without boundary by Lemma \ref{fiberprod}.

Now let $P = \mathcal{M}(B_1,B_2) \times_{B_2} \cdots \times_{B_n} 
\mathcal{M}(B_n,B')$ and assume for the purpose of induction that 
$\partial_-:P \rightarrow B_1$ is a submersion. This assumption
implies that the projection
$$
\pi_1:\mathcal{M}(B,B_1) \times_{B_1} P \rightarrow \mathcal{M}(B,B_1)
$$
is a submersion.  To see this, let $v \in T_x\mathcal{M}(B,B_1)$, 
and let $\alpha:(-\varepsilon,\varepsilon) \rightarrow 
\mathcal{M}(B,B_1)$ be a smooth curve such that $\alpha(0) = x$ 
and $\alpha'(0) = v$.  Since $\partial_-:P \rightarrow B_1$ is a 
submersion it is locally a projection map. Hence, if 
$\partial_+(\alpha(0)) \in B_1$ is contained in the image of 
$\partial_-:P \rightarrow B_1$, then for $t$ sufficiently small we 
have $\partial_+(\alpha(t)) \in \partial_-(P)$.  This implies that 
we can find a smooth curve $\gamma:(-\varepsilon,\varepsilon)
\rightarrow P$ such that $\partial_+(\alpha(t)) = \partial_-(\gamma(t))$
for $t$ sufficiently small. The smooth curve $\alpha \times \gamma$ 
then lies in $\mathcal{M}(B,B_1) \times_{B_1} P$ for $t$ sufficiently small
and satisfies $d\pi_1 (\alpha'(0) \times \gamma'(0)) = v$. 

To complete the proof of the corollary, simply note that
$\partial_-:\mathcal{M}(B,B_1) \rightarrow B$ is a submersion by Lemma 
\ref{beginsubmersioninterior} and hence the composition of this map
with $\pi_1$:
$$
\partial_-: \mathcal{M}(B,B_1) \times_{B_1} \mathcal{M}(B_1,B_2) \times_{B_2}
\cdots \times_{B_n} \mathcal{M}(B_n,B') \rightarrow B
$$
is also a submersion.

\begin{flushright}
$\Box$
\end{flushright}

\begin{lemma}\label{fiberedmanifold}
If $B$ and $B'$ are connected critical submanifolds and
$\sigma:P \rightarrow B$ is a smooth map from a compact
smooth manifold with corners $P$, then
$$
P \times_B \overline{\mathcal{M}}(B,B')
$$
is a compact smooth manifold with corners.
\end{lemma}

\smallskip\noindent
Proof:
Let $\Delta = \{(b,b) |\ b \in B \}$ be the diagonal in $B \times B$.
By Corollary \ref{beginsubmersion} the beginning point
map $\partial_-:\overline{\mathcal{M}}(B,B') \rightarrow B$
is a submersion and a stratum submersion.  Hence,
the map $\sigma \times \partial_-: P \times \overline{\mathcal{M}}(B,B')
\rightarrow B \times B$ is transverse and stratum transverse 
to the submanifold $\Delta \subset B \times B$. Thus,
$$
P \times_{B} \overline{\mathcal{M}}(B,B') = (\sigma \times \partial_-)^{-1}(\Delta)
$$
is a manifold with corners by Theorem \ref{preimagecorners}.
The space $P \times_{B} \overline{\mathcal{M}}(B,B')$ is compact
because it is a closed subset of the compact space
$P \times \overline{\mathcal{M}}(B,B')$.

\begin{flushright}
$\Box$
\end{flushright}

\smallskip\noindent
\textbf{Proof of Lemma \ref{basiscorners}:} The faces of $I^N$ are 
obviously compact smooth manifolds with corners, and for any two 
connected critical submanifolds $B$ and $B'$ the compactified moduli 
space $\overline{\mathcal{M}}(B,B')$ is a compact smooth manifolds 
with corners by Theorem \ref{compactification}. The rest of the elements 
in $C_p$ are constructed inductively from these elements, and hence 
Lemma \ref{fiberedmanifold} can be applied inductively to show that 
all the elements of $C_p$ are compact smooth manifolds with corners.
The orientations on these manifolds were defined previously.

\begin{flushright}
$\Box$
\end{flushright}


\section{Independence of the function}\label{continuationsection}

In this section we prove that the homology of the Morse-Bott chain complex
$(C_\ast(f),\partial)$ is independent of the Morse-Bott-Smale function
$f$, and hence isomorphic to the singular homology of $M$ with integer
coefficients. The method of proof follows standard arguments found in 
\cite{AusMor}, \cite{BarLag}, \cite{FloCoh}, \cite{FloAnI}, \cite{SchMor}, 
and \cite{WebThe}. The outline of the proof is as follows.
\nocite{FloSym} \nocite{FloMor}

Given two Morse-Bott-Smale functions $f_1,f_2:M \rightarrow \mathbb{R}$
we pick a smooth function $F_{21}: M \times \mathbb{R} \rightarrow \mathbb{R}$ 
meeting certain transversality requirements such that
\begin{eqnarray*}
\lim_{t \rightarrow -\infty} F_{21}(x,t) & = & f_1(x) + 1\\
\lim_{t \rightarrow +\infty} F_{21}(x,t) & = & f_2(x) - 1
\end{eqnarray*}
for all $x \in M$.  The compactified moduli spaces of gradient flow lines of $F_{21}$ 
(the \textit{time dependent} gradient flow lines) are used to define a 
chain map $(F_{21})_\Box:C_\ast(f_1) \rightarrow C_\ast(f_2)$, where 
$(C_\ast(f_k),\partial)$ is the Morse-Bott chain complex of $f_k$ for
$k=1,2$.

Next we consider the case where we have four Morse-Bott-Smale
functions $f_k:M \rightarrow \mathbb{R}$ where $k=1,2,3,4$, and
we pick a smooth function $H:M \times \mathbb{R} \times \mathbb{R}
\rightarrow \mathbb{R}$ meeting certain transversality requirements
such that
\begin{eqnarray*}
\lim_{s \rightarrow -\infty} \lim_{t \rightarrow -\infty} H(x,s,t) & = & f_1(x) + 2\\
\lim_{s \rightarrow +\infty} \lim_{t \rightarrow -\infty} H(x,s,t) & = & f_2(x)\\
\lim_{s \rightarrow -\infty} \lim_{t \rightarrow +\infty} H(x,s,t) & = & f_3(x)\\
\lim_{s \rightarrow +\infty} \lim_{t \rightarrow +\infty} H(x,s,t) & = & f_4(x) - 2\\
\end{eqnarray*}
for all $x \in M$. The compactified moduli spaces of gradient flow lines of $H$ 
are used to define a chain homotopy between $(F_{43})_\Box \circ 
(F_{31})_\Box$ and $(F_{42})_\Box \circ (F_{21})_\Box$ where
$(F_{lk})_\Box:C_\ast(f_k) \rightarrow C_\ast(f_l)$ is the map
defined above for $k,l=1,2,3,4$. In homology the map $(F_{kk})_\ast:
H_\ast(C_\ast(f_k),\partial) \rightarrow H_\ast(C_\ast(f_k),\partial)$
is the identity for all $k$, and hence 
\begin{eqnarray*}
(F_{12})_\ast \circ (F_{21})_\ast & = & (F_{11})_\ast \circ (F_{11})_\ast = id\\
(F_{21})_\ast \circ (F_{12})_\ast & = & (F_{22})_\ast \circ (F_{22})_\ast = id.
\end{eqnarray*}
Therefore,
$$
(F_{21})_\ast: H_\ast(C_\ast(f_1),\partial) \rightarrow H_\ast(C_\ast(f_2),\partial)
$$
is an isomorphism.


\subsection{Time dependent gradient flow lines}

Let $\rho:\mathbb{R} \rightarrow (-1,1)$ 
be a smooth strictly increasing function such that 
$\lim_{t \rightarrow -\infty} \rho(t) = - 1$ and
$\lim_{t \rightarrow +\infty} \rho(t)= 1$. The function
$\rho$ determines a function on the compactified real line
$\overline{\mathbb{R}} = \mathbb{R} \cup \{\pm \infty\}$ by defining
$\rho(-\infty) = -1$ and $\rho(+\infty) = 1$.  
(The compactified real line $\overline{\mathbb{R}}$ is a compact
smooth manifold with boundary with a single chart
$h:\overline{\mathbb{R}} \rightarrow [-1,1]$ given by
$h(t) = \frac{t}{\sqrt{1+t^2}}$, see Definition 2.1 of \cite{SchMor}.
Note that we could choose $\rho = h$.)

Let $f_1,f_2:M \rightarrow \mathbb{R}$ be Morse-Bott-Smale functions, and
let $B_i^{f_1}$ and $B_j^{f_2}$ denote the set of critical points of $f_1$
and $f_2$ respectively of index $i,j=0,\ldots ,m$. Let 
$F_{21}:M \times \mathbb{R} \rightarrow \mathbb{R}$ be a smooth 
function that is strictly decreasing in its second component
such that for some large $T \gg 0$ we have
$$
F_{21}(x,t) = \left\{ \begin{array}{llc}
f_1(x) - \rho(t) & \text{ if } & t < -T \\
h_t(x)   & \text{ if } & -T\leq t \leq T \\
f_2(x) - \rho(t) & \text{ if } & t > T
\end{array}\right.
$$
where $h_t(x)$ is an approximation to $\frac{1}{2}(T-t)(f_1(x) - \rho(t)) + 
\frac{1}{2}(T+t)(f_2(x) - \rho(t))$ that makes $F_{21}$ smooth.
(To simplify the notation in the following computations we will 
sometimes drop the subscript on $F_{12}$.)

Let $\varphi_\alpha:M \times \mathbb{R} \rightarrow M \times \mathbb{R}$ 
denote the flow associated to $-\nabla F_{21}$ with respect to the product
metric on $M \times \mathbb{R}$. Thus, $\varphi_\alpha(x,t) = 
(\gamma(\alpha), \sigma(\alpha)) \in M \times \mathbb{R}$ where
$$
(\gamma'(\alpha),\sigma'(\alpha)) = \left\{ \begin{array}{ll}
(-(\nabla f_1)(\gamma(\alpha)), \rho'(\alpha)) & \text{ if } t < -T \vspace{2 mm} \\
(-(\nabla f_2)(\gamma(\alpha)), \rho'(\alpha)) & \text{ if } t > T
\end{array}\right.
$$
and $(\gamma(0),\sigma(0)) = (x,t)$. The negative gradient flow
$\varphi_\alpha:M \times \mathbb{R} \rightarrow M \times \mathbb{R}$ of $F_{21}$
determines stable and unstable manifolds in $M \times \mathbb{R}$
of the critical submanifolds of $f_1$ and $f_2$:
\begin{eqnarray*}
W^s_F(B_j^{f_2}) & = & \{(x,t) \in M \times \mathbb{R} | \lim_{\alpha \rightarrow
+\infty} \pi_1(\varphi_\alpha(x,t)) \in B_j^{f_2} \}\\
W^u_F(B_i^{f_1}) & = & \{(x,t) \in M \times \mathbb{R} | \lim_{\alpha \rightarrow
-\infty} \pi_1(\varphi_\alpha(x,t)) \in B_i^{f_1} \}
\end{eqnarray*}
where $\pi_1: M \times \mathbb{R} \rightarrow M$ denotes projection onto the
first component. These stable and unstable manifolds inherit
orientations from the stable and unstable manifolds of $f_1$ and $f_2$, and
since both $f_1$ and $f_2$ satisfy the Morse-Bott-Smale transversality condition
we can choose the approximation $h_t(x)$ so that $W^u_F(B_i^{f_1}) \pitchfork
W^s_F(B_j^{f_2})$ is an oriented submanifold of $M \times \mathbb{R}$ for all
$i,j = 0, \ldots ,m$.

\begin{definition}
The moduli spaces of gradient flow lines of $F_{21}$ (the time dependent gradient
flow lines )are defined to be
$$
\mathcal{M}_F(B_i^{f_1},B_j^{f_2}) = (W^u_F(B_i^{f_1}) \cap W^s_F(B_j^{f_2}))/\mathbb{R}
$$
for all $i,j=0,\ldots ,m$ where the $\mathbb{R}$-action is defined by
the negative gradient flow.
\end{definition}

\begin{lemma}
For a generic function $F_{21}$ satisfying the conditions listed
above the moduli space $\mathcal{M}_F(B_i^{f_1},B_j^{f_2})$ is either empty
or an oriented smooth manifold of dimension $b_i^{f_1} + i-j$.
\end{lemma}

\smallskip\noindent
Proof:
The dimension of $W^u_F(B_i^{f_1})$ is $b_i^{f_1} + i + 1$ and the dimension of 
$W^s_F(B_j^{f_2})$ is $m+1-j$. Thus, the dimension of the
transverse intersection $W^u_F(B_i^{f_1}) \pitchfork W^s_F(B_j^{f_2})$ is 
$$
(b_i^{f_1} + i + 1) + (m+1-j) - (m+1) = b_i^{f_1} + i - j +1.
$$

\begin{flushright}
$\Box$
\end{flushright}

\noindent
Note: The preceding lemma implies that $\mathcal{M}_F(B_i^{f_1},B_j^{f_2}) = \emptyset$
whenever $b_i^{f_1} + i < j$.  This is a weaker statement than the one found
in Lemma \ref{weakly}, where the function is assumed to satisfy the Morse-Bott-Smale
transversality condition.

\smallskip
There are gluing and compactification results for time dependent gradient
flow lines analogous to Theorems \ref{gluing} and \ref{compactification}.
In order to distinguish between the time independent gradient flow lines and
the time dependent gradient flow lines, we will denote the moduli 
spaces of (time independent) gradient flow lines of $f_1$ by 
$\mathcal{M}_{f_1}(B_i^{f_1},B_j^{f_1})$.  Similarly, we will denote the 
(time independent) moduli spaces of gradient flow lines of $f_2$ by 
$\mathcal{M}_{f_2}(B_i^{f_2},B_j^{f_2})$. The proofs of the 
following two fundamental theorems follow from the results in 
Appendix \S A.2 and \S A.3 of \cite{AusMor}. See \cite{BouCom}, 
\cite{HofAGeI}, and \cite{HofAGeII} for an alternate approach to 
these theorems.

\begin{theorem}[Gluing]\label{timedepgluing}
For large $N \gg 0$  there are injective local diffeomorphisms
\begin{eqnarray*}
G_L:\mathcal{M}_{f_1}(B_i^{f_1},B_n^{f_1}) \times_{B_n^{f_1}} \mathcal{M}_F(B_n^{f_1},B_j^{f_2}) \times (N,\infty)
& \rightarrow & \mathcal{M}_F(B_i^{f_1},B_j^{f_2})\\
G_R:\mathcal{M}_F(B_i^{f_1},B_n^{f_2}) \times_{B_n^{f_2}} \mathcal{M}_{f_2}(B_n^{f_2},B_j^{f_2}) \times (N,\infty)
& \rightarrow & \mathcal{M}_F(B_i^{f_1},B_j^{f_2})
\end{eqnarray*}
onto the open ends of $\mathcal{M}_F(B_i^{f_1},B_j^{f_2})$.
\end{theorem}

\begin{theorem}[Compactification]\label{timedepcompactification}
The moduli space $\mathcal{M}_F(B_i^{f_1},B_j^{f_2})$ has a compactification 
$\overline{\mathcal{M}}_F(B_i^{f_1},B_j^{f_2})$, consisting of all the piecewise 
gradient flow lines from $B_i^{f_1}$ to $B_j^{f_2}$ (including both time 
dependent and time independent gradient flow lines), which is either 
empty or a compact oriented smooth manifold with corners of dimension 
$b_i^{f_1} + i-j$. The beginning and endpoint maps extend to smooth maps
on the compactified spaces.
\end{theorem}

The compactification of the moduli spaces of time dependent gradient
flow lines by piecewise gradient flow lines determines a degree
and boundary operator that satisfies the axioms for abstract
topological chains. To simplify the notation we will assume that
for each $i=0,\ldots ,m$ the components of $B_i^{f_1}$ are of the
same dimension. In general one needs to group the components
by their dimension and then define the degree and boundary
operator on each group.

\begin{definition}\label{timedepboundary}
The degree of $\overline{\mathcal{M}}_F(B_i^{f_1},B_j^{f_2})$ is defined
to be $b_i^{f_1} + i-j$, and the boundary operator 
$\partial \overline{\mathcal{M}}_F(B_i^{f_1},B_j^{f_2})$ is defined to
be $(-1)^{i+b_i^{f_1}}$ times the quantity
$$
\left( \sum_{n<i}\overline{\mathcal{M}}_{f_1}(B_i^{f_1},B_n^{f_1}) \times_{B_n^{f_1}} 
\overline{\mathcal{M}}_F(B_n^{f_1},B_j^{f_2}) \ - 
\sum_{j<n} \overline{\mathcal{M}}_F(B_i^{f_1},B_n^{f_2})
\times_{B_n^{f_2}} \overline{\mathcal{M}}_{f_2}(B_n^{f_2},B_j^{f_2}) \right)
$$
where $b_i^{f_1} = \mbox{dim }B_i^{f_1}$ and the fibered products are taken over 
the beginning and endpoint maps $\partial_-$ and $\partial_+$. 
If $B_n = \emptyset$, then the corresponding compactified moduli spaces
and fibered products are identified with zero. The boundary operator
extends to fibered products of compactified moduli spaces via Definition 
\ref{fiberedboundary}.
\end{definition}

\noindent
The proof of the following lemma is analogous to the proof of Lemma 
\ref{modulicomplex}.

\begin{lemma}
The degree and boundary operator for $\overline{\mathcal{M}}_F(B_i^{f_1},B_j^{f_2})$
satisfy the axioms for abstract topological chains, i.e. the boundary operator
is of degree $-1$ and it satisfies $\partial \circ \partial = 0$.
\end{lemma}


\subsection{Representing chain systems}

Intuitively, the chain map $(F_{21})_\Box:\tilde{C}_\ast(f_1) \rightarrow 
\tilde{C}_\ast(f_2)$ should be defined as follows. Given a singular 
$C_p^{f_1}$-space $\sigma:P \rightarrow B_i^{f_1}$ that is transverse and stratum
transverse to $\partial_-:\overline{\mathcal{M}}_F(B_i^{f_1},B_j^{f_2}) 
\rightarrow B_i^{f_1}$ for all $j=0,\ldots ,m$, we would like to define
$$
(F_{21})_\Box(\sigma) = \bigoplus_{j=0}^m \left(\partial_+\circ \pi_2:P \times_{B_i^{f_1}} 
\overline{\mathcal{M}}_F(B_i^{f_1},B_j^{f_2}) \rightarrow B_j^{f_2}\right) \in \bigoplus_{j=0}^m
S_{p+i-j}^\infty(B_j^{f_2})
$$
where $P \times_{B_i^{f_1}} \overline{\mathcal{M}}_F(B_i^{f_1},B_j^{f_2}) 
\stackrel{\pi_2}{\longrightarrow} \overline{\mathcal{M}}_F(B_i^{f_1},B_j^{f_2})
\stackrel{\partial_+}{\longrightarrow} B_j^{f_2}$ is defined as in 
Lemma \ref{partialj}. However, fibered products of the compactified 
moduli spaces of \textbf{time dependent} gradient flow lines with
singular cubes in the critical submanifolds of $f_1$ are not elements of the 
sets $\{C_p^{f_2}\}_{p \geq 0}$ of allowed domains for the singular topological 
chains in $\tilde{C}_\ast(f_2)$. So, to define the map rigorously 
we need to choose a topological chain in $S_{p+i-j}^\infty(B_j^{f_2})$ that 
represents the map $\partial_+ \circ \pi_2:P \times_{B_i^{f_1}} 
\overline{\mathcal{M}}_F(B_i^{f_1},B_j^{f_2}) \rightarrow B_j^{f_2}$. 
To do this we adapt the idea of representing chain systems from 
\cite{BarLag} to our setup.

For every $p \geq 0$ let $C_p^F$ be the set of $p$-dimensional components 
of the fibered products of the compactified moduli spaces of time dependent 
gradient flow lines of $F_{21}$ with elements of $\tilde{C}_\ast(f_1)$ that 
intersect the $\partial_-$ maps transversally and stratum transversally,
and let $S_p^F$ be the free abelian group generated by the elements of 
$C_p^F$, i.e. $S_p^F = \mathbb{Z}[C_p^F]$. Note that for any $p \geq 0$ 
the elements of $C_p^F$ are compact oriented smooth manifolds with corners, 
and hence they have fundamental classes in relative homology. Moreover, 
$(S_\ast^F,\partial)$ is a chain complex of abstract topological chains
by Lemma \ref{fiberedcomplex}.

\begin{definition}
A \textbf{representing chain system} for $\{C_p^F\}_{p \geq 0}$
is a family of smooth singular cubical chains $\{s_P \in S_p(P)|\ P \in
C_p^F\}_{p \geq 0}$ such that
\begin{enumerate}
\item The image of the smooth singular cubical chain $s_P$ in the relative
      cubical chain group $S_p(P,\partial P) = S_p(P)/S_p(\partial P)$
      (where $\partial P$ denotes the topological boundary of $P$) is a 
      cycle that represents the fundamental class of $P$.
\item If $\partial P = \sum_k n_k P_k \in S_{p-1}^F$ where $n_k = \pm 1$ and 
      $P_k \in C_{p-1}^F$, then
$$
\partial s_P = \sum_k n_k s_{P_k} \in S_{p-1}(P)
$$
where $\partial s_P$ denotes the boundary operator on smooth singular
cubical chains.
\end{enumerate}
\end{definition}

\begin{lemma}
There exists a representing chain system for $\{C_p^F\}_{p \geq 0}$.
\end{lemma}

\smallskip\noindent
Proof: The proof is by induction on the degree $p$.  For $p=0$ the
statement is obvious. So, assume that we have chosen representing
chain systems for all the elements in $C_j^F$ where $j < p$, and
let $P\in C_p^F$.  Consider the following diagram
$$
\xymatrix{
S_p^F \ar[r]^{\partial} \ar[d]^{R_p} & S_{p-1}^F \ar[d]^{R_{p-1}} \ar[r]^{\partial} & S_{p-2}^F \ar[d]^{R_{p-2}}\\
S_p(P) \ar[r]^{\partial} & S_{p-1}(P) \ar[r]^{\partial} & S_{p-2}(P)
}
$$
where $R_{p-1}$ and $R_{p-2}$ are defined by sending linear 
combinations of subsets of $P$ in $S_{p-1}^F$ and $S_{p-2}^F$ 
to their representing chains and all other elements to
zero. 

If $\partial P = \sum_k n_k P_k \in S_{p-1}^F$, then by the
induction hypothesis $R_{p-1}(\partial P) = \sum_k n_k s_{P_k}$ 
where $s_{P_k}\in S_{p-1}(P_k) \subseteq S_{p-1}(P)$ is
a smooth singular cubical chain whose image in
$S_{p-1}(P_k,\partial P_k)$ is a cycle that represents
the fundamental class of $P_k$.  Therefore, we can
chose an element $s_P \in S_p{P}$ such that $\partial s_P =
\sum_k n_k s_{P_k}$ and the image of $s_P$ in $S_p(P,\partial P)$
represents the fundamental class of $P$.  (See for instance
Lemma VI.9.1 of \cite{BreTop}.) Define $R_p(P) = s_P$.

\begin{flushright}
$\Box$
\end{flushright}


\subsection{The chain map $(F_{21})_\Box:C_\ast(f_1) \rightarrow C_\ast(f_2)$}

Assume that a representing chain system $\{s_P \in S_p(P)|\ P \in
C_p^F\}_{p \geq 0}$ has been chosen. The map $(F_{21})_\Box:
\tilde{C}_\ast(f_1)\rightarrow \tilde{C}_\ast(f_2)$ is defined 
as follows.

\begin{definition}\label{chainmapdef}
If $\sigma:P \rightarrow B_i^{f_1}$ is a singular $C_p^{f_1}$-space
in $B_i^{f_1}$ that intersects the beginning point map $\partial_-:
\overline{\mathcal{M}}_F(B_i^{f_1},B_j^{f_2}) \rightarrow B_i^{f_1}$
transversally and stratum transversally, then define
$$
(F_{21})_\Box(\sigma) = \bigoplus_{j=0}^m F_{21}^j(\sigma) \in \bigoplus_{j=0}^m
S_{p+i-j}^\infty(B_j^{f_2})
$$
where
$$
F_{21}^j(\sigma) = \sum_k n_k(\partial_+ \circ \pi_2 \circ s_{k})
$$
and
$$
\sum_k n_k s_{k} \in S_{p+i-j}(P \times_{B_i^{f_1}} \overline{\mathcal{M}}_F(B_i^{f_1},B_j^{f_2}))
$$
is the representing chain for $P \times_{B_i^{f_1}} \overline{\mathcal{M}}_F(B_i^{f_1},B_j^{f_2})$
if $P \times_{B_i^{f_1}} \overline{\mathcal{M}}_F(B_i^{f_1},B_j^{f_2}) \neq \emptyset$ and
$F_{21}^j(\sigma) = 0$ otherwise.
\end{definition}

\smallskip\noindent
Note: By general transversality arguments we may assume without loss of
generality that the map $\sigma:P \rightarrow B_i^{f_1}$ intersects the 
beginning point map $\partial_-:\overline{\mathcal{M}}_F(B_i^{f_1},B_j^{f_2})
\rightarrow B_i^{f_1}$ transversally and stratum transversally.
Also, to simplify the notation we will suppress the representing chain
systems and denote the chain $F_{21}^j(\sigma) \in S_{p+i-j}^\infty(B_j^{f_2})$
in the preceding definition by
$$
\partial_+ \circ \pi_2:P \times_{B_i^{f_1}} \overline{\mathcal{M}}_F(B_i^{f_1},B_j^{f_2}) 
\rightarrow B_j^{f_2}. 
$$

\noindent
The proof of the following lemma is similar to the proof of
Lemma \ref{boundarydegenerate}.

\begin{lemma}\label{inducedchainmaps}
The map $(F_{21})_\Box:\tilde{C}_\ast(f_1) \rightarrow \tilde{C}_\ast(f_2)$
induces a map $(F_{21})_\Box:C_\ast(f_1) \rightarrow C_\ast(f_2)$.
\end{lemma}

\begin{proposition}\label{chainmap}
The map $(F_{21})_\Box:C_\ast(f_1) \rightarrow C_\ast(f_2)$ is a chain 
map of degree zero. That is, the map preserves degree and
$\partial \circ (F_{21})_\Box = (F_{21})_\Box \circ \partial$.
\end{proposition}

\smallskip\noindent
Proof:
Let $\sigma:P \rightarrow B_i^{f_1}$ be a singular $C_p^{f_1}$-space in $B_i^{f_1}$ that
intersects the beginning point map $\partial_-:\overline{\mathcal{M}}_F(B_i^{f_1},B_j^{f_2})
\rightarrow B_i^{f_1}$ transversally and stratum transversally. The dimension 
of $P \times_{B_i^{f_1}} \overline{\mathcal{M}}_F(B_i^{f_1},B_j^{f_2})$ is $p+i-j$ 
by Lemma \ref{fiberprod} and Theorem \ref{timedepcompactification}, and thus 
$\partial_+ \circ \pi_2:P \times_{B_i^{f_1}} \overline{\mathcal{M}}_F(B_i^{f_1},B_j^{f_2}) 
\rightarrow B_j^{f_2}$ has Morse-Bott degree $p+i$ for all $j=0,\ldots ,m$.
Therefore, the map $(F_{21})_\Box:C_\ast(f_1) \rightarrow C_\ast(f_2)$ preserves degree.

We will now show that $(F_{21})_\Box(\partial(\sigma)) = \partial (F_{21})_\Box(\sigma)
\in \tilde{C}_\ast(f_2)$. Thus $(F_{21})_\Box(\partial(\sigma)) = \partial (F_{21})_\Box(\sigma)
\in C_\ast(f_2)$ by Lemma \ref{boundarydegenerate} and Lemma \ref{inducedchainmaps}. 
To simplify the notation we give the following computation in terms of abstract 
topological chains and assume that for every $i=0,\ldots ,m$ all the components 
of $B_i^{f_1}$ are of the same dimension $b_i^{f_1}$. In general the components need 
to be grouped by their dimensions and the computation repeated on each group.
The abstract topological chain associated to $\partial(\sigma) \in C_{p+i-1}(f_1)$ is
$$
(-1)^{p+i} \partial P + \sum_{k<i} P \times_{B_i^{f_1}} \overline{\mathcal{M}}_{f_1}(B_i^{f_1},B_k^{f_1}).
$$
Hence, by Definition \ref{chainmapdef} we have $(F_{21})_\Box(\partial(\sigma)) = $
$$
(-1)^{p+i} \sum_{j=0}^m \partial P \times_{B_i^{f_1}}
\overline{\mathcal{M}}_F(B_i^{f_1},B_j^{f_2}) + \sum_{j=0}^m \sum_{k<i} P \times_{B_i^{f_1}}
\overline{\mathcal{M}}_{f_1}(B_i^{f_1},B_k^{f_1}) \times_{B_k^{f_1}} \overline{\mathcal{M}}_F(B_k^{f_1},B_j^{f_2}).
$$
Using Definition \ref{chainmapdef} again we have $(F_{21})_\Box(\sigma) = \sum_{j=0}^m P \times_{B_i^{f_1}} 
\overline{\mathcal{M}}_F(B_i^{f_1},B_j^{f_2})$, and hence by Definition \ref{fiberedboundary}, 
Definition \ref{boundary}, and Lemma \ref{timedepboundary}
\begin{eqnarray*}
\partial_0 (F_{21})_\Box(\sigma) & = & (-1)^{p+i} \sum_{j=0}^m \partial(P \times_{B_i^{f_1}}
      \overline{\mathcal{M}}_F(B_i^{f_1},B_j^{f_2}))\\
& = & (-1)^{p+i} \sum_{j=0}^m \left[\partial P \times_{B_i^{f_1}}
      \overline{\mathcal{M}}_F(B_i^{f_1},B_j^{f_2}) + (-1)^{p+b_i^{f_1}} P \times_{B_i^{f_1}} \partial
      \overline{\mathcal{M}}_F(B_i^{f_1},B_j^{f_2})\right] \\
& = & (-1)^{p+i} \sum_{j=0}^m \partial P \times_{B_i^{f_1}}
      \overline{\mathcal{M}}_F(B_i^{f_1},B_j^{f_2}) + (-1)^{i+b_i^{f_1}} \sum_{j=0}^m 
      P \times_{B_i^{f_1}} \partial \overline{\mathcal{M}}_F(B_i^{f_1},B_j^{f_2}) \\
\end{eqnarray*}
\begin{eqnarray*}
& = & (-1)^{p+i} \sum_{j=0}^m \partial P \times_{B_i^{f_1}} \overline{\mathcal{M}}_F(B_i^{f_1},B_j^{f_2}) \\
&   & \qquad + \sum_{j=0}^m \sum_{k<i}  P \times_{B_i^{f_1}}
      \overline{\mathcal{M}}_{f_1}(B_i^{f_1},B_k^{f_1}) \times_{B_k^{f_1}} \overline{\mathcal{M}}_F(B_k^{f_1},B_j^{f_2}) \\
&   & \qquad - \sum_{j=0}^m \sum_{j<q} P \times_{B_i^{f_1}}
      \overline{\mathcal{M}}_F (B_i^{f_1},B_q^{f_2}) \times_{B_q^{f_2}} \overline{\mathcal{M}}_{f_2}(B_q^{f_2},B_j^{f_2}).
\end{eqnarray*}
Since
\begin{eqnarray*}
(\partial_1 \oplus \cdots \oplus \partial_m)((F_{21})_\Box(\sigma)) & = & 
\sum_{0 \leq k < j \leq m} P \times_{B_i^{f_1}} \overline{\mathcal{M}}_F(B_i^{f_1},B_j^{f_2}) 
\times_{B_j^{f_2}} \overline{\mathcal{M}}_{f_2}(B_j^{f_2},B_k^{f_2})\\
& = & \sum_{j=0}^m \sum_{j<q} P \times_{B_i^{f_1}} \overline{\mathcal{M}}_F 
      (B_i^{f_1},B_q^{f_2}) \times_{B_q^{f_2}} \overline{\mathcal{M}}_{f_2}(B_q^{f_2},B_j^{f_2})
\end{eqnarray*}
we have $(F_{21})_\Box(\partial(\sigma)) = (\partial_0 \oplus \cdots \oplus \partial_m) 
((F_{21})_\Box(\sigma)) \in \tilde{C}_\ast(f_2)$ and hence $(F_{21})_\Box(\partial(\sigma)) = \partial
(F_{21})_\Box(\sigma) \in C_\ast(f_2)$.

\begin{flushright}
$\Box$
\end{flushright}

\begin{corollary}\label{flowhomo}
The map $F_{21}:M \times \mathbb{R} \rightarrow \mathbb{R}$ induces a
homomorphism in homology
$$
(F_{21})_\ast:H_\ast(C_\ast(f_1),\partial) \rightarrow H_\ast(C_\ast(f_2),\partial) 
$$
which is independent of the choice of representing chain system.
\end{corollary}


\subsection{Chain homotopies from time dependent gradient flows}

\nocite{SalLec}

Assume now that we have four Morse-Bott-Smale
functions $f_k:M \rightarrow \mathbb{R}$, where $k=1,2,3,4$.
By Proposition \ref{chainmap} these four functions
determine chain maps
$$
(F_{lk})_\Box:C_\ast(f_k) \rightarrow C_\ast(f_l)
$$ 
for all $k,l=1,2,3,4$.  We will show that the moduli
spaces of gradient flow lines of a function
$H:M \times \mathbb{R} \times \mathbb{R} \rightarrow \mathbb{R}$
meeting certain transversality requirements can be used to construct a chain homotopy 
between $(F_{43})_\Box \circ (F_{31})_\Box$ and 
$(F_{42})_\Box \circ (F_{21})_\Box$.

Following \cite{AusMor}, \cite{FloSym}, \cite{SalLec}, and \cite{WebThe} we let 
$H:M \times \mathbb{R} \times \mathbb{R} \rightarrow \mathbb{R}$
be a smooth function that is strictly decreasing in its last two
components such that for some large $T \gg 0$ we have
$$
H(x,s,t) = \left\{ \begin{array}{lllll}
f_1(x) - \rho(s) - \rho(t) & \text{ if } & s < -T & \text{ and } & t < -T\\
f_2(x) - \rho(s) - \rho(t) & \text{ if } & s > T  & \text{ and } & t < -T\\
f_3(x) - \rho(s) - \rho(t) & \text{ if } & s < -T & \text{ and } & t > T\\
f_4(x) - \rho(s) - \rho(t) & \text{ if } & s > T  & \text{ and } & t > T
\end{array}\right.
$$
where $\rho:\mathbb{R} \rightarrow (-1,1)$ is a smooth strictly increasing 
function such that $\lim_{t \rightarrow -\infty} \rho(t) = - 1$ and
$\lim_{t \rightarrow \infty} \rho(t)= 1$. Taking the standard metric on
$\mathbb{R}$ and the product metric on $M \times \mathbb{R} \times \mathbb{R}$, 
the gradient flow of $H:M \times \mathbb{R} \times \mathbb{R}\rightarrow \mathbb{R}$
on the last two components of its domain can be pictured 
in $\overline{\mathbb{R}} \times \overline{\mathbb{R}}$ as follows.

\begin{center}
\includegraphics{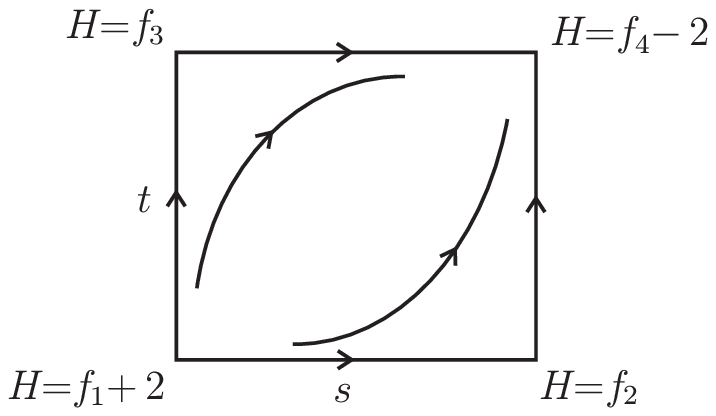}
\end{center}

Let $\varphi_\alpha:M \times \mathbb{R} \times \mathbb{R} \rightarrow 
M \times \mathbb{R} \times \mathbb{R}$ denote the flow associated to 
$-\nabla H$ with respect to the product metric on $M \times \mathbb{R}
\times \mathbb{R}$, and denote the stable and unstable manifolds
of the critical point sets $B_j^{f_4}$ and $B_i^{f_1}$ by
\begin{eqnarray*}
W^s_H(B_j^{f_4}) & = & \{(x,s,t) \in M \times \mathbb{R} \times \mathbb{R}|
\lim_{\alpha \rightarrow +\infty} \pi_1(\varphi_\alpha(x,s,t)) \in B_j^{f_4} \}\\
W^u_H(B_i^{f_1}) & = & \{(x,s,t) \in M \times \mathbb{R} \times \mathbb{R}|
\lim_{\alpha \rightarrow -\infty} \pi_1(\varphi_\alpha(x,s,t)) \in B_i^{f_1} \}
\end{eqnarray*}
where $\pi_1: M \times \mathbb{R} \times \mathbb{R} \rightarrow M$ denotes 
projection onto the first component. These stable and unstable manifolds inherit
orientations from the stable and unstable manifolds of $f_1$ and $f_4$, and
since both $f_1$ and $f_4$ satisfy the Morse-Bott-Smale transversality condition
we can choose $H$ so that $W^u_H(B_i^{f_1}) \pitchfork W^s_H(B_j^{f_4})$ 
is an oriented submanifold of $M \times \mathbb{R} \times \mathbb{R}$ 
for all $i,j = 0, \ldots ,m$.

The moduli spaces of gradient flow lines of $H$ are defined to be
$$
\mathcal{M}_H(B_i^{f_1},B_j^{f_4}) = (W^u_H(B_i^{f_1}) \cap W^s_F(B_j^{f_4}))/\mathbb{R}
$$
for all $i,j=0,\ldots ,m$ where the $\mathbb{R}$-action is defined by
the negative gradient flow, and we have
$$
\text{dim }\mathcal{M}_H(B_i^{f_1},B_j^{f_4}) = b_i^{f_1} + i - j + 1.
$$ 
Moreover, there are gluing and compactification results for the
moduli spaces of $H$ similar to those stated above. Hence, the moduli 
space $\mathcal{M}_H(B_i^{f_1},B_j^{f_4})$ has a compactification 
$\overline{\mathcal{M}}_H(B_i^{f_1},B_j^{f_4})$, consisting of all 
the piecewise gradient flow lines of $H$ from $B_i^{f_1}$ to 
$B_j^{f_4}$, which is a compact smooth manifold with corners. 
The beginning and endpoint maps $\partial_-$ and $\partial_+$ 
extend to smooth maps on the compactified space.

\begin{theorem}\label{chainhomotopy}
The function $H:M \times \mathbb{R} \times \mathbb{R} \rightarrow \mathbb{R}$
induces a chain homotopy
$$
(F_{43})_\Box \circ (F_{31})_\Box - (F_{42})_\Box \circ (F_{21})_\Box
= \partial H_\Box + H_\Box \partial.
$$
\end{theorem}

\smallskip\noindent
Proof:
Assume that the components of $B_i^{f_1}$ are of the same dimension
and that appropriate representing chains have been chosen.
The boundary of $\overline{\mathcal{M}}_H(B_i^{f_1},B_j^{f_4})$
is $(-1)^{i+b_i^{f_1}}$ times the quantity
\begin{eqnarray*}
&   & \sum_{k=0}^m \overline{\mathcal{M}}_{F_{21}}(B_i^{f_1}, B_k^{f_2}) \times_{B_k^{f_2}}
                   \overline{\mathcal{M}}_{F_{42}}(B_k^{f_2}, B_j^{f_4})\\
& - & \sum_{n=0}^m \overline{\mathcal{M}}_{F_{31}}(B_i^{f_1}, B_n^{f_3}) \times_{B_n^{f_3}}
                   \overline{\mathcal{M}}_{F_{43}}(B_n^{f_3}, B_j^{f_4})\\
& + & \sum_{n<i}   \overline{\mathcal{M}}_{f_1}(B_i^{f_1}, B_n^{f_1}) \times_{B_n^{f_1}}
                   \overline{\mathcal{M}}_{H}(B_n^{f_1}, B_j^{f_4})\\
& + & \sum_{j<n}   \overline{\mathcal{M}}_{H}(B_i^{f_1}, B_n^{f_4}) \times_{B_n^{f_4}}
                   \overline{\mathcal{M}}_{f_4}(B_n^{f_4}, B_j^{f_4}).\\
\end{eqnarray*}
If $\sigma:P \rightarrow B_i^{f_1}$ satisfies the required transversality
requirements, then $H_\Box(\sigma)\in C_\ast(f_4)$ is represented by the
abstract topological chain
$$
\sum_{j=0}^m P \times_{B_i^{f_1}} \overline{\mathcal{M}}_H(B_i^{f_1},B_j^{f_4})
$$
and $\partial H_\Box(\sigma)\in C_\ast(f_4)$ is represented by the
abstract topological chain
$$
\sum_{j=0}^m \left( \partial_0 \left(P \times_{B_i^{f_1}} \overline{\mathcal{M}}_H(B_i^{f_1},B_j^{f_4})\right) +
\sum_{k<j} P \times_{B_i^{f_1}} \overline{\mathcal{M}}_H(B_i^{f_1},B_j^{f_4})
      \times_{B_j^{f_4}} \overline{\mathcal{M}}_{f_4}(B_j^{f_4},B_k^{f_4})\right).
$$
For the $\partial_0$ terms in the above sum we have
$\partial_0 \left(P \times_{B_i^{f_1}} \overline{\mathcal{M}}_H(B_i^{f_1},B_j^{f_4})\right) = $
$$
(-1)^{p+i+1} \left( \partial P \times_{B_i^{f_1}} \overline{\mathcal{M}}_H(B_i^{f_1},B_j^{f_4}) +
(-1)^{p+b_i^{f_1}} P \times_{B_i^{f_1}} \partial \overline{\mathcal{M}}_H(B_i^{f_1},B_j^{f_4}) \right)
$$
where $\partial\overline{\mathcal{M}}_H(B_i^{f_1},B_j^{f_4})$ is the sum given above.
In addition, the term $H_\Box(\partial \sigma)$ is represented by the abstract topological 
chain
$$
\sum_{j=0}^m \left( (-1)^{p+i}\partial P \times_{B_i^{f_1}} \overline{\mathcal{M}}_H(B_i^{f_1},B_j^{f_4}) +
\sum_{n<i} P \times_{B_i^{f_1}} \overline{\mathcal{M}}_{f_1}(B_i^{f_1},B_n^{f_1}) \times_{B_n^{f_1}} 
\overline{\mathcal{M}}_H(B_n^{f_1},B_j^{f_4}) \right).
$$
Hence, $\partial H_\Box(\sigma) + H_\Box(\partial \sigma)$ is represented by the
sum from $j=0$ to $m$ of the abstract topological chains
$$\hspace{-0.7cm}
\sum_{n=0}^m P \times_{B_i^{f_1}} \overline{\mathcal{M}}_{F_{31}}(B_i^{f_1}, B_n^{f_3}) 
\times_{B_n^{f_3}} \overline{\mathcal{M}}_{F_{43}}(B_n^{f_3}, B_j^{f_4}) -
\sum_{k=0}^m P \times_{B_i^{f_1}} \overline{\mathcal{M}}_{F_{21}}(B_i^{f_1}, B_k^{f_2})
\times_{B_k^{f_2}} \overline{\mathcal{M}}_{F_{42}}(B_k^{f_2}, B_j^{f_4}).
$$

\begin{flushright}
$\Box$
\end{flushright}

\begin{corollary}\label{identity}
For every Morse-Bott-Smale function $f:M \rightarrow \mathbb{R}$
there exists a smooth function 
$F:M \times \mathbb{R} \rightarrow \mathbb{R}$ such that
$$
F_\ast:H_\ast(C_\ast(f),\partial) \rightarrow H_\ast(C_\ast(f),\partial)
$$
is the identity homomorphism.  
\end{corollary}

\smallskip\noindent
Proof:
Let $F:M \times \mathbb{R} \rightarrow \mathbb{R}$ be
defined as $F(x,t) = f(x) - \rho(t)$ for all $(x,t) \in
M \times \mathbb{R}$. If $\sigma:P \rightarrow B_i^{f}$ is a
singular $C_p^{f}$-space in $B_i^{f}$, then
$$
F_\Box(\sigma) = \sum_{j \leq i} \left( \partial_+:
P \times_{B_i^{f}} \overline{\mathcal{M}}_{F}(B_i^{f},B_j^{f})
\rightarrow B_j^{f} \right)
$$
since the Morse-Bott index of a Morse-Bott-Smale function is
decreasing along its gradient flow lines by Lemma \ref{weakly}.
The moduli space $\overline{\mathcal{M}}_{F}(B_i^{f},B_i^{f})$ is 
diffeomorphic to $B_i^{f}$, and the beginning point map
$
\partial_-:\overline{\mathcal{M}}_{F}(B_i^{f},B_i^{f}) \rightarrow B_i^{f}
$
is the identity map with respect to this diffeomorphism.  Hence,
$$
P \times_{B_i^{f}} \overline{\mathcal{M}}_{F}(B_i^{f},B_i^{f}) \approx
\{(x,y) \in P \times B_i^{f} |\ \sigma(x) = y\} \approx P,
$$
and the endpoint map $\partial_+:P \times_{B_i^{f}} 
\overline{\mathcal{M}}_{F}(B_i^{f},B_i^{f}) \rightarrow B_i^{f}$
agrees with $\sigma:P \rightarrow B_i^{f}$.  Therefore,
$$
F_\Box(\sigma) = \sigma + \sum_{j < i} \left( \partial_+:
P \times_{B_i^{f}} \overline{\mathcal{M}}_{F}(B_i^{f},B_j^{f})
\rightarrow B_j^{f} \right).
$$
Note that if $f$ is a Morse-Smale function or a constant function, 
then the terms in the above sum with $j<i$ are all zero in $C_\ast(f)$,
and thus the map $F_\Box:C_\ast(f) \rightarrow C_\ast(f)$
is the identity at the level of chains.

In general, let $\tilde{F}:M \times \mathbb{R} \rightarrow \mathbb{R}$
be a smooth function that is strictly decreasing in its second 
component such that for some large $T \gg 0$ we have
$$
\tilde{F}(x,t) = \left\{ \begin{array}{llc}
f(x) - \rho(t) & \text{ if } & t < -\frac{1}{2}T \\
h_t(x)   & \text{ if } & -\frac{1}{2}T\leq t \leq \frac{1}{2}T \\
f(x) - \rho(t) & \text{ if } & t > \frac{1}{2}T
\end{array}\right.
$$
where the approximation $h_t(x)$ is chosen so that the 
non-constant time dependent gradient flow lines of $\tilde{F}$ all flow 
into $W^s(B_0^{f}) \times \mathbb{R}$. This is possible because 
$\text{dim }W^s(B_0^{f}) = m$ and every piecewise gradient flow
line of $f$ can be extended to a piecewise gradient flow line that 
ends in $B_0^{f}$. If $\sigma:P \rightarrow B_i^{f}$ is a singular 
$C_p^{f}$-space in $B_i^{f}$, then
$$
\tilde{F}_\Box(\sigma) = \sigma + \left(\partial_+: P \times_{B_i^{f}} 
\overline{\mathcal{M}}_{\tilde{F}}(B_i^{f},B_0^{f}) \rightarrow B_0^{f} \right).
$$
Now, let $f_1 = f_2 = f_4 = f$, and let $f_3\equiv 0$ be a constant 
function.  Choose 
$$
F_{31}(x,t)  = \left\{ \begin{array}{llc}
f(x) - \rho(t) & \text{ if } & t < -T \\
\mu(t)f(x) - \rho(t) & \text{ if } & -T\leq t \leq -\frac{1}{2}T \\
0 - \rho(t) & \text{ if } & t > -\frac{1}{2}T
\end{array}\right.
$$
where $\mu:\mathbb{R} \rightarrow [0,1]$ is a smooth decreasing 
function that is $1$ for $t \leq -T$ and $0$ for $t \geq -\frac{1}{2}T$. 
If we choose 
$$
F_{43}(x,t)  = \left\{ \begin{array}{llc}
0 - \rho(t) & \text{ if } & t < -T \\
(1-\mu(t))f(x) - \rho(t) & \text{ if } & -T\leq t \leq -\frac{1}{2}T \\
h_t(x)   & \text{ if } & -\frac{1}{2}T \leq t \leq \frac{1}{2}T \\
f(x) - \rho(t) & \text{ if } & t > \frac{1}{2}T
\end{array}\right.
$$
then for all $i=0,\ldots ,m$ we have diffeomorphisms
$$
\overline{\mathcal{M}}_{F_{31}}(B_i^{f},M) \times_{M}
\overline{\mathcal{M}}_{F_{43}}(M,B_0^{f}) \approx
\overline{\mathcal{M}}_{\tilde{F}}(B_i^{f},B_0^{f}),
$$
and hence, $(F_{43})_\Box \circ (F_{31})_\Box = \tilde{F}_\Box$.
Thus, if we choose $F_{21} = F_{42} = \tilde{F}$, then 
Theorem \ref{chainhomotopy} implies that
$
\tilde{F}_\ast = \tilde{F}_\ast \circ \tilde{F}_\ast.
$

If $\sum_k n_k\sigma_k$ is a cycle in $C_\ast(f)$, then we have
$$
\tilde{F}_\ast \left(\left[\sum_k n_k\sigma_k\right]\right) = \left[\sum_k n_k\sigma_k\right] + [C]
$$
where $C$ is a sum of singular spaces in $B_0^f$ of the
form $\partial_+: P_k \times_{B_i^{f}} \overline{\mathcal{M}}_{\tilde{F}}(B_{i_k}^{f},B_0^{f})
\rightarrow B_0^f$.
Hence, 
$$
(\tilde{F}_\ast \circ \tilde{F}_\ast) \left(\left[\sum_k n_k\sigma_k\right]\right) = 
\left[\sum_k n_k\sigma_k\right] + 2[C],
$$
and the fact that $\tilde{F}_\ast = \tilde{F}_\ast \circ \tilde{F}_\ast$
implies that $[C] = 0$.  Therefore, 
$$
\tilde{F}_\ast \left(\left[\sum_k n_k\sigma_k\right]\right) = \left[\sum_k n_k\sigma_k\right],
$$
and $\tilde{F}_\ast:H_\ast(C_\ast(f),\partial) \rightarrow H_\ast(C_\ast(f),\partial)$
is the identity homomorphism.  

\begin{flushright}
$\Box$
\end{flushright}

\begin{corollary}\label{canonical}
For any two Morse-Bott-Smale functions $f_1,f_2:M \rightarrow \mathbb{R}$
the time-dependent gradient flow lines from $f_1$ to $f_2$ determine
a canonical homomorphism 
$$
(F_{21})_\ast:H_\ast(C_\ast(f_1),\partial) \rightarrow H_\ast(C_\ast(f_2),\partial), 
$$
i.e. the map $(F_{21})_\ast$ is independent of the choice of the function
$F_{21}:M \times \mathbb{R} \rightarrow \mathbb{R}$ and the representing
chain system used to define the chain map $(F_{21})_\Box:C_\ast(f_1)
\rightarrow C_\ast(f_2)$.
\end{corollary}

\smallskip\noindent
Proof:
Let $f_2 = f_3 = f_4$ in Theorem \ref{chainhomotopy}, let 
$F_{21}:M \times \mathbb{R} \rightarrow \mathbb{R}$ and
$F_{31}:M \times \mathbb{R} \rightarrow \mathbb{R}$ be two
functions that define time-dependent gradient flow
lines from $f_1$ to $f_2 = f_3$, and let $F_{22}$ and
$F_{33}$ be the functions from Corollary \ref{identity}.
Theorem \ref{chainhomotopy} implies
$$
(F_{33})_\ast \circ (F_{31})_\ast = (F_{22})_\ast \circ (F_{21})_\ast
$$
and since $(F_{22})_\ast = (F_{33})_\ast = id$, we have
$(F_{31})_\ast = (F_{21})_\ast$.

\begin{flushright}
$\Box$
\end{flushright}

\begin{corollary}\label{functorial}
For any four Morse-Bott-Smale functions $f_k:M \rightarrow \mathbb{R}$, where
$k=1,2,3,4$, the canonical homomorphisms satisfy
$$
(F_{43})_\ast \circ (F_{31})_\ast = (F_{42})_\ast \circ (F_{21})_\ast
$$
and
$$
(F_{32})_\ast \circ (F_{21})_\ast =  (F_{31})_\ast.
$$
\end{corollary}

\smallskip\noindent
Proof:  This follows immediately from Theorem \ref{chainhomotopy} and Corollary
\ref{identity}.

\begin{flushright}
$\Box$
\end{flushright}

\begin{theorem}\label{homologyindependence}
The homology of the Morse-Bott chain complex $(C_\ast(f),\partial)$ is
independent of the Morse-Bott-Smale function $f:M \rightarrow \mathbb{R}$.
Therefore, 
$$
H_\ast(C_\ast(f),\partial) \approx H_\ast(M;\mathbb{Z}).
$$
\end{theorem}

\smallskip\noindent
Proof:  Let $f_1,f_2:M \rightarrow \mathbb{R}$ be Morse-Bott-Smale 
functions.  If we take $f_3 = f_1$ in Corollary \ref{functorial},
then we see that $(F_{12})_\ast \circ (F_{21})_\ast =  (F_{11})_\ast = id$
by Corollaries \ref{identity} and \ref{canonical}. Similarly,
$(F_{21})_\ast \circ (F_{12})_\ast = id$.

\begin{flushright}
$\Box$
\end{flushright}

\bibliographystyle{amsxport}
\bibliography{books,papers}

\end{document}